\documentclass{IEEEtran}

\newcommand{\mr}{\mathrm}                   % Roman, non-italic
\newcommand{\veg}[1]{\bm{#1}}               % Bold for geometrical / physical vectors
\newcommand{\mat}[1]{\mathsfbfit{#1}}       % NxN dimensional matrix
\newcommand{\mLam}{\mat \Lambda}            % Loop matrix
\newcommand{\mSig}{\mat \Sigma}             % Star matrix
\newcommand{\T}{\mr{T}}

\usepackage[utf8]{inputenc}

\usepackage[T1]{fontenc}
\usepackage{trimclip}
\usepackage{float}
\usepackage{tikz}
\usepackage{pgfplots}
\usepackage{pgfplotstable}
\usepackage{siunitx}
\pgfplotsset{compat=newest}
\pgfplotsset{every mark/.append style={solid}}
\pgfplotsset{every axis/.append style={
             label style={font=\small},
             tick label style={font=\small}}}
\usepgfplotslibrary{units}

\usepackage{kbordermatrix}
\renewcommand{\kbldelim}{(}
\renewcommand{\kbrdelim}{)}

\usepackage{newtxtext}
\usepackage{amsthm}
\DeclareMathAlphabet{\mathbfsf}{\encodingdefault}{\sfdefault}{bx}{n}
\usepackage{bm}
\usepackage[OMLmathsfit]{isomath}

\usepackage{envmath}
\usepackage{mathtools}
\usepackage{amssymb}
\usepackage{esint}

\usepackage{multirow}
\usepackage{hyperref}
\usepackage{cleveref}
\allowdisplaybreaks[1]

\usepackage{comment}
\usepackage{makecell}

\usepackage{cite}
\usepackage[]{todonotes}

\usepgfplotslibrary{external}
\makeatletter
\renewcommand{\todo}[2][]{\tikzexternaldisable\@todo[#1]{#2}\tikzexternalenable}
\makeatother

\newcommand{\blue}[1]{#1}

\usepackage{amsmath}
\usepackage{subfig}

\renewcommand{\vec}[1]{\mathsfbfit{#1}}

\def\BibTeX{{\rm B\kern-.05em{\sc i\kern-.025em b}\kern-.08em T\kern-.1667em\lower.7ex\hbox{E}\kern-.125emX}}

\def\PLH{\mat{P}^{\boldsymbol{\Lambda H}}}
\def\PS{\mat{P}^{\boldsymbol{\Sigma}}}
\def\PSH{\mathbb{P}^{\boldsymbol{\Sigma H}}}
\def\PL{\mathbb{P}^{\boldsymbol{\Lambda}}}

\def\updelta{\delta}
\def\upsigma{\sigma}

\begin{document}
\title{Low-Frequency Stabilizations of the PMCHWT Equation for Dielectric and Conductive Media: On a Full-Wave Alternative to Eddy-Current Solvers}
\author{Viviana Giunzioni, \IEEEmembership{Graduate Student Member, IEEE}, Alberto Scazzola, Adrien Merlini, \IEEEmembership{Senior Member, IEEE}, and Francesco P. Andriulli, \IEEEmembership{Fellow, IEEE}

\thanks{
%This paragraph of the first footnote will contain the date on  which you submitted your paper for review. 
This work was supported by the European Innovation Council (EIC) through the European Union’s Horizon Europe research Programme under Grant 101046748 (Project CEREBRO) and by the European Union – Next Generation EU within the PNRR project ``Multiscale modeling and Engineering Applications'' of the Italian National Center for HPC, Big Data and Quantum Computing (Spoke 6) – PNRR M4C2, Investimento 1.4 - Avviso n. 3138 del 16/12/2021 - CN00000013 National Centre for HPC, Big Data and Quantum Computing (HPC) - CUP E13C22000990001.}
\thanks{Viviana Giunzioni and Francesco P. Andriulli are with the Department of Electronics and Telecommunications, Politecnico di Torino, 10129 Turin, Italy (e-mail: francesco.andriulli@polito.it).}
\thanks{Adrien Merlini is with the Microwave Department, IMT Atlantique, 29238 Brest, France (e-mail: adrien.merlini@imt-atlantique.fr).}
%\thanks{T. C. Author is with  the Electrical Engineering Department, University of Colorado, Boulder, CO  80309 USA, on leave from the National Research Institute for Metals,  Tsukuba, Japan (e-mail: author@nrim.go.jp).}
}

\maketitle

\begin{abstract}
We propose here a novel stabilization strategy for the PMCHWT equation that cures its frequency and conductivity related instabilities and is obtained by leveraging quasi-Helmholtz projectors. The resulting formulation is well-conditioned in the entire low-frequency regime, including the eddy current one, and can be applied to arbitrarily penetrable materials, ranging from dielectric to conductive ones. In addition, by choosing the rescaling coefficients of the quasi-Helmholtz components appropriately, we prevent the typical loss of accuracy occurring at low frequency in the presence of inductive and capacitive type magnetic frill excitations, commonly used in circuit modeling to impose a potential difference.
Finally, leveraging on quasi-Helmholtz projectors instead than on the standard Loop-Star decomposition,
our formulation is also compatible with most fast solvers and is amenable to multiply connected geometries, without any computational overhead for the search for the global loops of the structure.
The efficacy of the proposed preconditioning scheme when applied to both simply and multiply connected geometries is corroborated by numerical examples.
\end{abstract}

\begin{IEEEkeywords}
%Enter key words or phrases in alphabetical  order, separated by commas. For a list of suggested keywords, send a blank e-mail to keywords@ieee.org or visit \underline {http://www.ieee.org/organizations/pubs/ani\_prod/keywrd98.txt}
Eddy currents, preconditioning, full-wave, multiply
connected, quasi-Helmholtz
decomposition.
\end{IEEEkeywords}

\section{Introduction}
\label{sec:introduction}
\IEEEPARstart{T}{he} electromagnetic modeling of dielectric penetrable objects is crucial in a variety of applications, ranging from antenna design to micro- and nano-electronic applications \cite{balanis2012advanced, dirks1996quasistationary}. 
The boundary element method (BEM) allows to recast this problem in a set of integral equations \cite{sauter2011boundary}. Differently from other popular techniques, its use automatically enforces radiation conditions and only requires the discretization of surfaces separating different materials \cite{sauter2011boundary,steinbach2008numerical}. These characteristics make the BEM one of the most favorable and widespread tools to numerically solve the electromagnetic scattering problem.

The PMCHWT (Poggio-Miller-Chang-Harrington-Wu-Tsai) equation \cite{poggio1973integral} and the M\"{u}ller equation \cite{muller1958grundprobleme} are among the most common integral formulations for penetrable bodies, and are respectively of the first- and second-kind. In spite of its favorable stability properties, the M\"{u}ller equation in its standard, non-conforming, discretization has been shown to yield less accurate results than the PMCHWT formulation, especially in high-contrast scenarios \cite{bogaert2014lowfrequency,yla-oijala2008analysis,vantwout2022boundary}. Unfortunately, because of its first-kind nature, the PMCHWT formulation suffers from a dense-discretization ill-conditioning \cite{adrian2021electromagnetic,cools2011Calderon}. In addition, the PMCHWT equation is also plagued by a severe low-frequency breakdown \cite{adrian2021electromagnetic,beghein2017lowfrequency}, coupled to further instabilities related to the conductivity characterizing the obstacle medium.

To circumvent these issues, \textit{ad-hoc} solvers relying on quasi-static approximations of the Maxwell's system are usually employed to address the low-frequency regime. In particular, the simulation of eddy currents (i.e., conductive currents generated by a time-varying magnetic field \cite{hiptmair2007boundary,kriezis1992eddy}), of great interest for many industrial applications \cite{garcia-martin2011nondestructive}, has historically required the use of quasi-static solvers that neglect displacement currents \cite{rucker1995various}. One of the most severe drawback related to these approximations is the challenge of handling the simulation of multi-scale scenarios which restricts their use in real case, complex apparatuses.

%In recent years, the search for a full-wave formulation, well-conditioned and stable toward statics, as accurate as the PMCHWT equation, has been topic of deep research.
Several preconditioning techniques based on the Calder\'{o}n identities have been proposed for the PMCHWT formulation, relying on the fact that the PMCHWT operator, when properly discretized, is a valid preconditioner for itself \cite{niino2012Calderon, suyan2010comparative,cools2011Calderon}. While immune from the dense-discretization breakdown, these formulations still suffer from contrast dependent instabilities \cite{cools2011Calderon} and undesired current cancellations at very low frequency, which dramatically affect the accuracy of the simulation outcomes \cite{beghein2017lowfrequency,chhim2020eddy}. An appropriate rescaling of the quasi-Helmholtz decomposed solenoidal and quasi-irrotational components of the system has been proven effective in addressing this issue \cite{adrian2021electromagnetic}. However, when the quasi-Helmholtz decomposition is performed by means of the Loop-Star change of basis
%because of the differential nature of the Loop and Star to RWG (Rao-Wilton-Glisson) functions transformation matrices
the condition number behavior with respect to the mesh refinement is even worsened \cite{adrian2021electromagnetic}. Also, the computationally costly explicit detection of the global-loop quasi-harmonic subspace is required \cite{adrian2021electromagnetic}.

Alternatively, the quasi-Helmholtz components of the system can be retrieved by projection, that is, by multiplication of the operator matrix with the quasi-Helmholtz projectors, introduced in \cite{andriulli2012loopstar,andriulli2013wellconditioned} and proved effective at curing the PMCHWT equation's low frequency breakdown in \cite{chhim2020eddy}. Their use brings along two significant advantages with respect to the standard Loop-Star decomposition, which are (I) the quasi-harmonic subspace is included in the solenoidal subspace spanned by the solenoidal projector, so that the global-loops detection is not required, and (II) the use of the quasi-Helmholtz projector does not affect the condition number behavior against discretization \cite{adrian2021electromagnetic}. 

By leveraging on this tool, we develop here a new preconditioning strategy capable of stabilizing the PMCHWT equation with respect to the physical parameters of the problem in the whole low-frequency regime, including the eddy-current one \cite{hiptmair2007boundary}, and, contextually, preventing detrimental current cancellations in presence of both inductive and capacitive type magnetic frill excitations. Built upon the standard PMCHWT equation, our full-wave formulation allows for seamless transitions between low and high frequencies, as well as between dielectric and conductive scatterers, resulting in a versatile and computationally efficient solver, compatible with state-of-the-art acceleration strategies \cite{hackbuschsparse, hackbusch2000sparse}.

The paper is organized as follows: after setting the necessary background and notation in \Cref{sec:background}, Sections \ref{sec:cond_analysis} and \ref{sec:dominantcomponent} aim at analyzing the shortcomings of the PMCHWT formulation. In particular, Section \ref{sec:cond_analysis}, presents an analysis of the asymptotic scalings of the quasi-Helmholtz components of the system as a function of frequency and of the conductivity of the scattering body. Section \ref{sec:dominantcomponent}, is focused on the dominant component analysis of the equivalent currents induced by a voltage generator modelled as a magnetic frill, which is common type of excitation employed in circuit simulation.  
Then, in \Cref{sec:precondForm}, we define a quasi-Helmholtz projectors based preconditioning strategy and determine the multiplicative coefficients used in the rescaling of the quasi-Helmholtz components of the system, allowing for the stabilization of the PMCHWT equation with respect to both the electrical length and the conductivity of the scatterer and for the correct recovery of the currents dominant components. The numerical results in \Cref{sec:num_res} aim at illustrating the good conditioning properties of the proposed formulation applied to both simply and multiply connected geometries, with respect to both frequency and conductivity and the accuracy of the results. Finally, concluding remarks will be given in \Cref{sec:conclusion}.

\section{Background and Notation}
\label{sec:background}
Let $\Omega_1\subset \mathbb{R}^3$ be a closed, bounded conductive region with boundary $\Gamma = \partial \Omega_1$ characterized by the outward pointing normal $\hat{\veg n}$, surrounded by the complementary vacuum region $\Omega_0 = \mathbb{R}^3\backslash\bar{\Omega_1}$, with permittivity $\epsilon_0$ and permeability $\mu_0$. The domain $\Omega_1$, of characteristic length $L$, is filled with a linear, homogeneous, and isotropic conductive medium of permeability $\mu_1 = \mu_r\mu_0$ and complex permittivity $\epsilon_1 = \epsilon_r'\epsilon_0-j\sigma/\omega$, where $\sigma$ is the conductivity, and $\omega$ is the angular frequency. The skin depth parameter, or penetration length of the field inside the scatterer, is function of frequency and conductivity as $\updelta\coloneqq\sqrt{2/(\omega\sigma\mu_1)}$.
%In the following, the region $\Omega_1$ can be simply or multiply connected.

Given a wave number $k$ and the three-dimensional free-space Green's function associated to it \cite{steinbach2008numerical}, $G_k(\veg r,\veg r') \coloneqq \mathrm{e}^{-\mathrm{j}kR}/(4\pi R)$, with $R \coloneqq \|\veg r-\veg r'\|$, the electric and magnetic field integral operators (EFIO and MFIO), $\mathcal{T}_k$ and $\mathcal{K}_k$, are defined as \cite{steinbach2008numerical, sauter2011boundary}
\begin{align}
(\mathcal{T}_k\,\boldsymbol{f})(\veg r) &\coloneqq -\mathrm{j}k (\mathcal{T}_{A,k}\,\boldsymbol{f})(\veg r) + \frac{1}{\mathrm{j}k} (\mathcal{T}_{\Phi,k}\,\boldsymbol{f})(\veg r)\,,\\
(\mathcal{K}_k\,\boldsymbol{f})(\veg r) &\coloneqq \hat{\veg n}\times \text{p.v.}\int_\Gamma \nabla G_k(\veg r,\veg r') \times\boldsymbol{f}(\veg r') \mathrm{d} \veg r'\,,\\(\mathcal{T}_{A,k}\,\boldsymbol{f})(\veg r) &\coloneqq \hat{\veg n}\times \int_\Gamma G_k(\veg r,\veg r') \boldsymbol{f}(\veg r') \mathrm{d} \veg r' \,,\\
(\mathcal{T}_{\Phi,k}\,\boldsymbol{f})(\veg r) &\coloneqq \hat{\veg n}\times \nabla \int_\Gamma G_k(\veg r,\veg r') \nabla'\cdot\boldsymbol{f}(\veg r') \mathrm{d} \veg r'\,,
\end{align}
where p.v. indicates Cauchy principal value.
When an incident electromagnetic field $(\veg E^i,\veg H^i)$ impinges on $\Omega_1$, the resulting electric and magnetic surface current densities $\veg j_s \coloneqq \hat{\veg n}\times \veg H$ and $\veg m_s \coloneqq -\hat{\veg n}\times \veg E$ satisfy the PMCHWT equation \cite{poggio1973integral},
\begin{equation}
    \begin{pmatrix}
        \eta_0 \mathcal{T}_{k_0}+\eta_1 \mathcal{T}_{k_1} & - (\mathcal{K}_{k_0}+\mathcal{K}_{k_1}) \\
        (\mathcal{K}_{k_0}+\mathcal{K}_{k_1}) &  \frac{1}{\eta_0} \mathcal{T}_{k_0}+\frac{1}{\eta_1} \mathcal{T}_{k_1}
    \end{pmatrix} \begin{pmatrix}
        \veg j_s \\ \veg m_s
    \end{pmatrix} = \begin{pmatrix}
        -\hat{\veg n}\times \veg E^i \\ -\hat{\veg n}\times \veg H^i
    \end{pmatrix}\,,
    \label{eqn:PMCHWT}
\end{equation}
where $k_{0,1}\coloneqq\omega\sqrt{\epsilon_{0,1}\mu_{0,1}}$ is the wave number in $\Omega_{0,1}$ and $\eta_{0,1}\coloneqq\sqrt{\mu_{0,1}/\epsilon_{0,1}}$ is the characteristic impedance of the exterior or interior medium.

To discretize and numerically solve \eqref{eqn:PMCHWT}, we approximate the surface $\Gamma$ with a mesh of planar triangular elements with average edge length $h$, over which we can define different sets of basis functions. In the following, we will consider Rao-Wilton-Glisson (RWG) basis functions. For their definition and visual representation, we refer to \cite{adrian2021electromagnetic} and references therein.
By expanding the unknown current densities as linear combinations of RWG functions, as
\begin{equation}
    \veg j_s = \sum_{n=1}^{N_e}j_n \veg f_n, \quad \veg m_s = \sum_{n=1}^{N_e}m_n \veg f_n\,,
\end{equation}
where $N_e$ is the number of edges in the mesh, and by testing the resulting equations with curl-conforming rotated RWG functions, following a conforming Petrov-Galerkin discretization procedure, we obtain the linear system of equations
\begin{equation}
\begin{pmatrix}
\mat T_u & -\mat K \\
\mat K &\mat T_l
\end{pmatrix}
\begin{pmatrix}
\vec j \\
\vec m
\end{pmatrix}=
\begin{pmatrix}
\vec e \\
\vec h
\end{pmatrix}\,,\quad \quad \begin{pmatrix}
\mat T_u & -\mat K \\
\mat K &\mat T_l
\end{pmatrix}  \eqcolon \mat Z
\label{eqn:PMCHWT_discr}
\end{equation}
where
\begin{align}
    \mat K &\coloneqq \mat K_{k_0}+\mat K_{k_1}\,,\\
    \mat T_u &\coloneqq \eta_0 \mat T_{k_0}+\eta_1 \mat T_{k_1}\,,\\
    \mat T_l &\coloneqq \frac{1}{\eta_0} \mat T_{k_0}+ \frac{1}{\eta_1} \mat T_{k_1}\,.
\end{align}
The elements of the above matrices discretizing the EFIO and MFIO read
\begin{align}
(\mat T_k)_{mn} &\coloneqq -\mathrm{j}k(\mat T_{A,k})_{mn}+\frac{1}{\mathrm{j}k}(\mat T_{\Phi,k})_{mn}\,,\\
(\mat K_k)_{mn} &\coloneqq \left( \hat{\veg n}\times \veg f_m, \mathcal{K}_k(\veg f_n) \right)_{L^2(\Gamma)}\,,\\
(\mat T_{A,k})_{mn} &\coloneqq \left( \hat{\veg n}\times \veg f_m, \mathcal{T}_{A,k}(\veg f_n) \right)_{L^2(\Gamma)}\,,\\
(\mat T_{\Phi,k})_{mn} &\coloneqq \left( \hat{\veg n}\times \veg f_m, \mathcal{T}_{\Phi,k}(\veg f_n) \right)_{L^2(\Gamma)}\,.
\end{align}
The right-hand-side (RHS) resulting from the testing procedure is
\begin{align}
    (\vec e)_m &\coloneqq \left( \hat{\veg n}\times \veg f_m, - \hat{\veg n}\times \veg E^i \right)_{L^2(\Gamma)}\,,\\  (\vec h)_m &\coloneqq \left( \hat{\veg n}\times \veg f_m, - \hat{\veg n}\times \veg H^i \right)_{L^2(\Gamma)}\,,
\end{align}
while the elements of the unknown vector are simply given by $(\vec j)_n \coloneqq j_n$ and $(\vec m)_n \coloneqq m_n$.
In the following we will focus on a rescaled version of the equation, that is, on matrix $\bar{\mat Z}$, which reads
\begin{equation}
\underbracket{
\begin{pmatrix}
\mat T_u/\eta_0 & -\mat K \\
\mat K &\eta_0\mat T_l
\end{pmatrix}}_{\bar{\mat Z}}
\begin{pmatrix}
\sqrt{\eta_0}\vec j \\
\vec m / \sqrt{\eta_0}
\end{pmatrix}=
\begin{pmatrix}
\vec e / \sqrt{\eta_0} \\
\sqrt{\eta_0} \vec h
\end{pmatrix}\,,
\label{eqn:PMCHWT_discr_diag}
\end{equation}
%\todo[inline]{AM: Why are you calling things diag precond? Are you using the actual operator daigonal? RIGHT}
where
\begin{align}
    \mat T_u/\eta_0 &= \mat T_{k_0} -\mathrm{j}k_0\mu_r\mat T_{A,k_1}+\frac{1}{\mathrm{j}k_0\epsilon_r}\mat T_{\Phi,k_1}\,,\\
    \eta_0\mat T_l &= \mat T_{k_0} -\mathrm{j}k_0\epsilon_r\mat T_{A,k_1}+\frac{1}{\mathrm{j}k_0\mu_r}\mat T_{\Phi,k_1}\,.
\end{align}
This rescaling leads to an effective conductivity stabilization at mid-frequencies, as will be shown in the following, but cannot prevent the dramatic increase of the conditioning towards lower frequencies.

We introduce at this point the mixed Gram matrices $\mat G$ and $\mathbb{G} = -\mat G^\T$, defined as $(\mat G)_{mn} \coloneqq \left( \hat{\veg n}\times \veg f_m, \veg g_n\right)_{L^2(\Gamma)}$ and $(\mathbb{G})_{mn} \coloneqq \left( \hat{\veg n}\times \veg g_m, \veg f_n\right)_{L^2(\Gamma)}$,
where $g_n$ is the dual with respect to the RWG function defined on the $n$-th edge. In this work, we employed Buffa-Christianses (BC) functions but other possibilities are feasible (see \cite{adrian2021electromagnetic} and references therein for the definition of BC functions and for further details).

%In the remaining part of this Section, we will introduce the elements required for the stabilization of the PMCHWT equation. 
In \Cref{sec:cond_analysis} we will apply the Loop-Star decomposition technique to the PMCHWT operator for analyzing its spectral properties, and we denote with $\mat \Lambda$, $\mat H$, and $\mat \Sigma$ the transformation matrices from the loop, global loop, and star subspaces to the RWG subspace. For their explicit definition and interpretation, we refer the reader to \cite{adrian2021electromagnetic, andriulli2012loopstar}. We also define the block matrix $\mat A$ as $\mat A \coloneqq \left( \mat\Lambda \quad \mat H \quad \mat \Sigma \right)$, that can be used to perform the required decomposition. It is worth noticing that we are only introducing matrix $\mat H$ for the purpose of the analysis, but its explicit (potentially costly) evaluation will not be required for the implementation of the proposed formulation.
The main building blocks of the actual preconditioning scheme presented in this work are the quasi-Helmholtz projectors, for both primal and dual spaces, defined as
\begin{align}
\PS &\coloneqq \mat \Sigma (\mat \Sigma^\T \mat \Sigma)^+ \mat \Sigma^\T\,,\quad\ \PLH \coloneqq \mat I -\PS\,,\\
\PL &\coloneqq \mat \Lambda (\mat \Lambda^\T \mat \Lambda)^+ \mat \Lambda^\T\,,\quad\quad \PSH \coloneqq \mat I -\PL\,,
\end{align}
where the subscript $^+$ indicates the Moore-Penrose pseudoinvserse. In particular, $\PS$ and $\PL$ are the non-solenoidal projectors for primal and dual functions; by complementarity, $\PLH$ and $\PSH$ are the solenoidal projectors for primal and dual functions, that means that they are projectors onto the solenoidal subspace---which includes the quasi-harmonic subspace.
By leveraging algebraic multigrid preconditioning techniques, these projectors can be applied to a vector in quasi-linear complexity \cite{adrian2021electromagnetic}, characteristic which makes them compatible with standard fast solvers.

\section{Conditioning Analysis of the PMCHWT Equation}
\label{sec:cond_analysis}
This section aims at studying some of the spectral properties of the PMCHWT system matrix $\bar{\mat Z}$ in the low-frequency regime and for different conductivity levels $\sigma$ of the body $\Omega$. The analysis is based on the determination of the asymptotic scalings towards low frequency of the quasi-Helmholtz components of the decomposed system, where the decomposition is performed by means of loop and star transformation matrices, for analytical purposes only.

First, we introduce the scalar parameters
\begin{align}
    \chi &\coloneqq k_0 L\,,\\
    \gamma &\coloneqq \sqrt{\omega\epsilon_0/\sigma}\,,\\
    \xi &\coloneqq \sqrt{2/\mu_r}L/\updelta\,,
\end{align}
such that $\chi = \gamma \xi$.
Before moving to the study of the interior Green's function kernel, strictly dependent on the frequency-conductivity regime considered (the nomenclature of which has been borrowed from \cite{dirks1996quasistationary}), we start by analyzing the asymptotic behavior of matrices $\mat T_{A,k_0}$, $\mat T_{\Phi,k_0}$, and $\mat K_{k_0}$ for vanishing frequency, as well as their loop-star components. 
In the low-frequency limit, for $\chi\rightarrow 0$, the exterior Green's function and its gradient can be expanded as
\begin{align}
     G_{k_0}(\veg r,\veg r') &= \frac{1}{4\pi R}\left[\ 1 - \mathrm{j}\frac{R}{L}\chi - \frac{R^2}{2 L^2}\chi^2 + \mathcal{O}(\chi^3) \right]\,,\label{eqn:Gexp}\\
     %\mathrm{j}\frac{R^3}{6 L^3}\chi^3 +\mathcal{O}(\chi^4) \right]\,,\label{eqn:Gexp}\\
     \nabla G_{k_0}(\veg r,\veg r') &= \frac{1}{4\pi}\left[\ \nabla\left(\frac{1}{R}\right) - \frac{\chi^2}{2 L^2}\nabla R + \mathcal{O}(\chi^3) \right]\,.\label{eqn:GradGexp}%\mathrm{j}\frac{\chi^3}{6 L^3}\nabla R^2 + \mathcal{O}(\chi^4) \right]\,,\label{eqn:GradGexp}
\end{align}
From these expansions, it can be recognized that, for $\chi\rightarrow 0$,
\begin{align}
\|\mat T_{A,k_0}/L\| &= \mathcal{O}(1)\,,\quad \|\mat K_{k_0}\| = \mathcal{O}(1)\,,\\ \|\mat T_{\Phi,k_0}L\| &= \mathcal{O}(1)\,,\quad \|\mat K_{\text{ext},k_0}\| = \mathcal{O}\left(\chi^2\right)
\,,
\end{align}
where $\mat K_{\text{ext},k} = \mat K_{k} - \mat K_{0}$. They result in the Loop-Star decompositions
\begin{align}
    \mat A^\T \left(\mat T_{k_0}\right) \mat A &= \mathcal{O}\begin{pmatrix}
    \chi & \chi & \chi \\
    \chi & \chi & \chi \\
    \chi & \chi & \chi^{-1}
    \end{pmatrix}\,,\\
    \mat A^\T \left(\mat K_{k_0}\right) \mat A &= \mathcal{O}\begin{pmatrix}
    \chi^2 & \chi^2 & 1 \\
    \chi^2 & 1 & 1 \\
    1 & 1 & 1
    \end{pmatrix}\,,
    \label{eqn:k0LoopStar}
\end{align}
where we have enforced the cancellation of the scalar potential contribution when tested against or applied to a solenoidal function \cite{andriulli2013wellconditioned}
\begin{equation}
    \mat \Lambda^\T \mat T_{\Phi,k} = \boldsymbol{0}\,,\quad \mat H^\T \mat T_{\Phi,k} = \boldsymbol{0}\,,\quad \mat T_{\Phi,k}\mat \Lambda= \boldsymbol{0}\,,\quad \mat T_{\Phi,k}\mat H= \boldsymbol{0} \label{eqn:Tpcanc}
\end{equation}
and the cancellation of the static magnetic field operator $\mat K_{0}$ \cite{bogaert2011low}
\begin{equation}
    \mat \Lambda^\T \mat K_{0} \mat \Lambda = \boldsymbol{0}\,,\quad \mat \Lambda^\T \mat K_{0} \mat H = \boldsymbol{0}\,,\quad \mat H^\T \mat K_{0} \mat \Lambda = \boldsymbol{0}\,. \label{eqn:K0canc}
\end{equation}
The decompositions in \eqref{eqn:k0LoopStar} will be useful to determine the overall Loop-Star decomposed PMCHWT matrix $\mat Z_{\mat \Lambda \mat H \mat \Sigma}$,
\begin{equation}
    \mat Z_{\mat \Lambda \mat H \mat \Sigma} \coloneqq \begin{pmatrix}
        \mat A^\T & \boldsymbol{0} \\ \boldsymbol{0} & \mat A^\T
    \end{pmatrix}\bar{\mat Z}\begin{pmatrix}
        \mat A  & \boldsymbol{0} \\ \boldsymbol{0} & \mat A 
    \end{pmatrix}\,.
\end{equation}

We proceed now to the definition of the three low-frequency regimes considered in this work and to the interior kernel analysis.

\subsection{Analysis in the Quasi-Static Regime}
\label{sec:analysisQSR}
The quasi-static regime is characterized by slowly varying fields, for which the quasi-stationary condition $\chi \ll 1$ holds true. In addition, the ratio between the angular frequency and the conductivity should be high enough to make any conductive current inside the object negligible with respect to the displacement current induced by the impinging field. In formulae, the quasi-static regime, abbreviated to QSR, is characterized by
\begin{equation}
    \begin{cases}
     \chi \ll 1\\
    \gamma \gg 1
    \end{cases}\,.
\end{equation}

The asymptotic behavior of the blocks of the Loop-Star decomposed PMCHWT equation, allowing the study of the stability properties of the formulation, can be deduced from the low-frequency asymptotic expansion of the interior Green's function and of its gradient, under the assumption $\gamma \gg 1$. This corresponds to analyzing the low-frequency limit along a constant$-\gamma$ line (blue arrow in \Cref{fig:regimes}) with $\gamma \gg 1$.

Under these conditions both the real and the imaginary part of the interior wave number vanish. In particular, writing the expression of the interior wavenumber in the form
\begin{equation}
    k_1 = k_0\sqrt{\epsilon_r'\mu_r} \sqrt{1-\mathrm{j}\frac{\sigma}{\omega\epsilon_0\epsilon_r'}}\,,
\end{equation}
shows that, for $\chi \rightarrow 0$, $\Re\{k_1 L\}=\mathcal{O}(\chi)$ and $\Im\{k_1 L\}=\mathcal{O}(\chi/\gamma^2)$.
This fact allows us to leverage on the small argument Taylor expansions of the Green's function and of its gradient to determine the low-frequency scalings
\begin{align}
\|\mat T_{A,k_1}/L\| &= \mathcal{O}(1)\,,\quad \|\mat K_{k_1}\| = \mathcal{O}(1)\,,\\ \|\mat T_{\Phi,k_1}L\| &= \mathcal{O}(1)\,,\quad \|\mat K_{\text{ext},k_1}\| = \mathcal{O}\left(\chi^2\right)\,,
\label{eqn:QSR_opscaling}
\end{align}
leading to the quasi-Helmholtz decomposition
\begin{equation}
    \mat Z_{\mat \Lambda \mat H \mat \Sigma} = \mathcal{O}\begin{pmatrix}
        \chi & \chi & \chi & \chi^2 & \chi^2 & 1 \\
        \chi & \chi & \chi & \chi^2 & 1 & 1 \\
        \chi & \chi & \chi^{-1} & 1 & 1 & 1 \\
        \chi^2 & \chi^2 & 1 & \chi & \chi & \chi \\
        \chi^2 & 1 & 1 & \chi & \chi & \chi \\
        1 & 1 & 1 & \chi & \chi & \chi^{-1}
    \end{pmatrix}\,.
    \label{eqn:LSPMCHWT_QS}
\end{equation}
By applying the Gershgorin circle theorem, we recognize that the condition number of the PMCHWT equation grows when decreasing the frequency, resulting in the so-called low-frequency breakdown of the formulation. Moreover, we infer that, as in the EFIE case \cite{andriulli2013wellconditioned}, the pathological frequency behavior of the condition number is in part due to the opposed frequency scaling of the vector and scalar potential components of the electric operator $\mathcal{T}_{k}$. Numerical results indicate that the condition number of the PMCHWT formulation increases as $\omega^{-2}$ in the low-frequency limit, as for the EFIE.

\subsection{Analysis in the Eddy-Current-Free Eddy-Current Regime}
The eddy-current regime is a quasi-stationary regime in which the typical time-scale is longer than the electric relaxation time constant $\tau = \epsilon_0 / \sigma$, i.e. $\omega\tau \ll 1$, equivalent to $\gamma\ll1$.
In addition, the eddy-current-free subset of the eddy-current regime is characterized by a further condition on the skin depth parameter, which should be much larger with respect to the characteristic size of the object, resulting in the condition $\xi\ll1$ in the regime under study.
Therefore, in the eddy-current-free eddy-current regime, abbreviated as ECFR, the frequency, geometrical and electrical parameters of the scatterer are such that
\begin{equation}
    \begin{cases}
     \chi \ll 1\\
     \gamma \ll 1
     \\\xi \ll 1
    \end{cases}\,.
\end{equation}

In the eddy current regime ($ \chi \ll 1$ and $\gamma \ll 1$) the interior wave number
\begin{equation}
    k_1 = \frac{1-\mathrm{j}}{\updelta}\sqrt{1+\mathrm{j}\frac{\omega\epsilon_0\epsilon_r'}{\sigma}} 
\end{equation}
can be approximated as $k_1\simeq(1-\mathrm{j})/\updelta$. Hence, under the further assumption $\xi \ll 1$, both the real and the imaginary parts of $k_1 L$ vanish as $\xi$.

In this regime, we want to analyze the low-frequency limit of the PMCHWT equation by keeping the ratio $\xi/\gamma$ constant (i.e., constant conductivity), corresponding to the black arrow in \Cref{fig:regimes}. In this limit, $\Re\{k_1 L\}=\mathcal{O}(\chi/\gamma)$ and $\Im\{k_1 L\}=\mathcal{O}(\chi/\gamma)$ for $\chi \rightarrow 0$.
Once again we can take advantage of the small argument expansions of the Green's function and of its gradient to determine the scalings 
\begin{align}
\|\mat T_{A,k_1}/L\| &= \mathcal{O}(1)\,,\quad \|\mat K_{k_1}\| = \mathcal{O}(1)\,,\\ \|\mat T_{\Phi,k_1}L\| &= \mathcal{O}(1)\,,\quad \|\mat K_{\text{ext},k_1}\| = \mathcal{O}\left(\chi/\gamma\right)^2\,.
\end{align}

The Loop-Star decomposed PMCHWT equation in the low-$\chi$ limit follows as
\begin{equation}
    \mat Z_{\mat \Lambda \mat H \mat \Sigma} = \mathcal{O}\begin{pmatrix}
        \chi & \chi & \chi & \chi^2\gamma^{-2} & \chi^2\gamma^{-2} & 1 \\
        \chi & \chi & \chi & \chi^2\gamma^{-2} & 1 & 1 \\
        \chi & \chi & \chi^{-1} & 1 & 1 & 1 \\
        \chi^2\gamma^{-2} & \chi^2\gamma^{-2} & 1 & \chi\gamma^{-2} & \chi\gamma^{-2} & \chi\gamma^{-2} \\
        \chi^2\gamma^{-2} & 1 & 1 & \chi\gamma^{-2} & \chi\gamma^{-2} & \chi\gamma^{-2} \\
        1 & 1 & 1 & \chi\gamma^{-2} & \chi\gamma^{-2} & \chi^{-1}
    \end{pmatrix}\,.
    \label{eqn:LSPMCHWT_ECFR}
\end{equation}
As in the previous case, the application of the Gershgorin circle theorem shows a low-frequency breakdown of the equation, as the condition number grows at least as $\omega^{-1}$ in the low-frequency limit.

Indeed, the matrix eigenvalues lie in three Gershgorin circles with radius $\mathcal{O}(1)$ and centered in values which are $\mathcal{O}(\chi)$, $\mathcal{O}(\chi^{-1})$, and $\mathcal{O}(\chi\gamma^{-2})$. Also notice then, when fixing the frequency, $\chi\gamma^{-2}=\sigma\eta_0L$ ranges from $\chi$ to $\chi^{-1}$ when the conductivity varies from $\upsigma_{\text{low}}\coloneqq\omega\epsilon_0$ (corresponding to $\gamma = 1$) to $\upsigma_{\text{high}}\coloneqq(\omega L^2\mu_1)^{-1}$ (corresponding to $\xi = 1$).

\subsection{Analysis in the Skin-Effect-Dominated Eddy-Current Regime}
The skin-effect-dominated eddy-current-regime, abbreviated as SEDR, is characterized by the relations
\begin{equation}
    \begin{cases}
     \chi \ll 1\\
     \gamma \ll 1\\
     \xi \gg 1
    \end{cases}\,.
\end{equation}
In this case, we aim at studying the asymptotic beahviour of the equation in the low-frequency limit along a constant$-\xi$ line, with  $\xi \gg 1$, corresponding to the red arrow in \Cref{fig:regimes}.

The eddy-current conditions ($ \chi \ll 1$ and $\gamma \ll 1$) guarantee the validity of the approximation of the interior wave number as $k_1 \simeq (1-\mathrm{j})/\updelta$ mentioned before. The fundamental difference with respect to the other low-frequency regimes analyzed above is that in this case the product $k_1L$ does not vanish in the low-frequency limit, preventing the use of the Maclaurin series of the Green's function and its gradient for the evaluation of the asymptotic behavior of the operators. However, by noticing that the Green's function behavior is dominated by the exponential decay $\mathrm{e}^{-R/\updelta}$, we can infer the asymptotic scalings
\begin{align}
\|\mat T_{A,k_1}/L\| &= \xi^{-1}\mathcal{O}(1)\,,\quad \|\mat K_{k_1}\| = \xi^{-1}\mathcal{O}(1)\,,\\ \|\mat T_{\Phi,k_1}L\| &= \xi^{-1}\mathcal{O}(1)\,,\quad \|\mat K_{\text{ext},k_1}\| = \mathcal{O}(1)\,
\end{align}
for $\chi\rightarrow 0$, where $\xi^{-1}=\gamma/\chi$ is constant.

Indeed, as far as the vector and scalar electric potentials are concerned, it can be seen that, in the low-skin-depth limit, only the self and near interactions give rise to a non negligible contribution as the skin depth becomes smaller than the characteristic size of the mesh triangles $h$. Hence, we can grasp an intuition of the above scalings by noticing that the Green's function kernel in modulo can be approximated as $|G_k(\veg r,\veg r')|\simeq \mathrm{e}^{-R/\updelta}/(4 \pi R)$ and that the integral of such function over an infinite plane, which is a good estimate of the integral over the triangle $T$ due to the exponential decay of the function for source points inside $T$ and testing points outside $T$, is
\begin{equation}
    \int_{\mathbb{R}^2}\frac{\mathrm{e}^{-|\veg r|/\updelta}}{4\pi|\veg r|\,}\mathrm{d} \veg r =\updelta/2 \,,
\end{equation}
determining the scalings of the norms $\|\mat T_{A,k_1}\|$ and $\|\mat T_{\Phi,k_1}\|$ with the skin depth. Also, in the same low-skin-depth-limit, the magnetic field integral operator matrix becomes negligible, as the kernel exactly cancels for testing points lying on the source triangle. Finally, each of $\mat T_{A,k_1}$, $\mat T_{\Phi,k_1}$, and $\mat K_{k_1}$, denoted here with the placeholder $\mat M$, can be written as the sum of a matrix accounting for the interactions between basis functions defined on non overlapping domains, $\mat M_{\text{far}}$, and the remaining $\mat M_{\text{near}}$, that is $\mat M = \mat M_{\text{near}} + \mat M_{\text{far}}$. It can be shown that, in the low-frequency limit, $\|\mat M_{\text{far}}\| = e^{-p\xi}\mathcal{O}(1)$, where $pL > 0$ represents the minimum distance between two non-overlapping domains of basis functions. Hence, $\|\mat M\|\le \|\mat M_{\text{near}}\|+\|\mat M_{\text{far}}\|=\xi^{-1}\mathcal{O}(1) + e^{-p\xi}\mathcal{O}(1)$.

The accurate evaluation of the inner integrals for the definition of the operator matrices is not trivial and may require the implementation of different integration techniques specifically tailored to highly lossy media. In our work, we employed the integration scheme presented in \cite{zhiguoqian2007generalized} to produce accurate numerical results.

Due to the condition $\gamma \ll 1$ the diagonal blocks of $\bar{\mat Z}$ can be approximated as
\begin{align}
    \mat T_u/\eta_0 &\simeq \mat T_{k_0} -\mathrm{j}k_0\mu_r\mat T_{A,k_1}+\frac{1}{\sigma\eta_0 L}(\mat T_{\Phi,k_1}L)\\
    \mat T_l\cdot\eta_0 &\simeq \mat T_{k_0} -\sigma\eta_0L(\mat T_{A,k_1}/L)+\frac{1}{\mathrm{j}k_0\mu_r}\mat T_{\Phi,k_1}\,,
\end{align}
It holds $\sigma\eta_0 L = \chi\gamma^{-2} \gg \chi^{-1}$ (following from the condition $\xi \gg 1$). As a result, the asymptotic behavior of the upper diagonal block is determined by $\mat T_{k_0}$, while the behavior of the lower diagonal block is determined by the dominating term $\sigma\eta_0L(\mat T_{A,k_1}/L)$, for which it holds
\begin{equation}
    \|\sigma\eta_0L(\mat T_{A,k_1}/L)\|=\mathcal{O}({\chi\gamma^{-2}\xi^{-1}})=\mathcal{O}({\gamma^{-1}})\,.
\end{equation}

Under these conditions, the Loop-Star decomposition of the PMCHWT equation is
\begin{equation}
    \mat Z_{\mat \Lambda \mat H \mat \Sigma} = \mathcal{O}\begin{pmatrix}
        \chi & \chi & \chi & 1 & 1 & 1 \\
        \chi & \chi & \chi & 1 & 1 & 1 \\
        \chi & \chi & \chi^{-1} & 1 & 1 & 1 \\
        1 & 1 & 1 & \gamma^{-1} & \gamma^{-1} & \gamma^{-1} \\
        1 & 1 & 1 & \gamma^{-1} & \gamma^{-1} & \gamma^{-1} \\
        1 & 1 & 1 & \gamma^{-1} & \gamma^{-1} & \gamma^{-1}
    \end{pmatrix}\,.
    \label{eqn:SEDR_loopstarPMCHWT}
\end{equation}
Also in the SEDR the formulation is plagued by the low-frequency breakdown, that is, the condition number of $\bar{\mat Z}$ increases at least as $\omega^{-1}$ in the low-frequency limit.
Moreover, a further instability with respect to $\sigma$ can be recognized in the high-conductivity limit, related to the branch of singular values accumulating at $\mathcal{O}(\gamma^{-1})$, with $\gamma^{-1}$ always higher that $\chi^{-1}$ in the SEDR, resulting in a condition number increase when fixing the frequency and increasing the conductivity of the scatterer.

Before moving to the definition of a quasi-Helmholtz projectors based preconditioning strategy capable of curing the instabilities identified so far, we first propose a dominant current components analysis, allowing to determine which are the current components required for the accurate evaluation of the electric and magnetic fields inside and outside the scatterer in the low-frequency limit. This preliminary analysis will indeed result in further requirements to be satisfied in the design of the preconditioning strategy.

\section{Dominant Component Analysis}
\label{sec:dominantcomponent}
While a stable and low matrix condition number is required for the convergence of an iterative solver \cite{quarteroni2008numerical}, it is not sufficient to guarantee the accuracy of the solution \cite{adrian2021electromagnetic,andriulli2013wellconditioned,chhim2020eddy}.  
Indeed, numerical losses of significant digits at low frequencies can provoke the corruption of some quasi-Helmholtz current components, which, in turn, can lower the accuracy of the scattered fields. In this section, we present a current dominant component analysis aimed at determining which are the current components required for the correct reconstruction of the electric and magnetic fields inside and outside the body $\Omega$. This information will constitute a further constraint in the design of the preconditioning strategy, which will have to allow for the correct recovery of these required current components.

\begin{figure}[!t]
\centerline{\includegraphics[width=0.7\columnwidth]{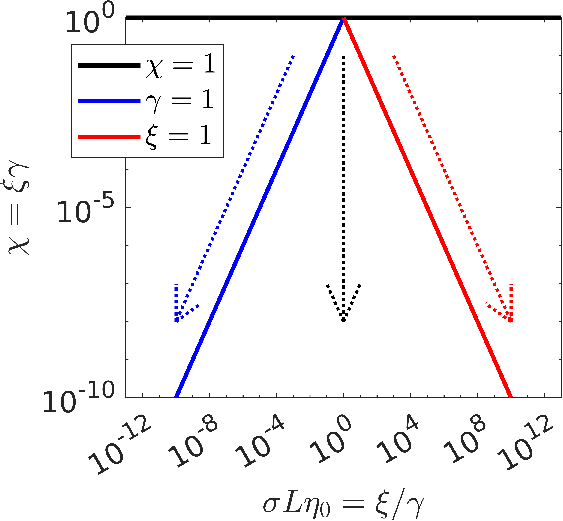}}
\caption{Graphical representation of the low-frequency limit in the three regimes: the blue arrow represents the limit direction in the QSR, the black arrow in the ECFR, and the red arrow in the SEDR.}
\label{fig:regimes}
\end{figure}

The analysis will leverage on the low-frequency asymptotic expansion of the PMCHWT equation along the low-frequency limits previously identified, which are, along a constant$-\gamma$ line with $\gamma \gg 1$ for the QSR, along a constant$-(\xi/\gamma)$  with $\xi \gg 1$ and $\gamma \gg 1$ for the ECFR, and along a constant$-\xi$ line with $\xi \gg 1$ for the SEDR.  These three limits are graphically represented as arrows in \Cref{fig:regimes}.

The dominant component analysis will require the determination of the asymptotic scalings of the loop, star and quasi-harmonic components of the system, both in their real and imaginary parts, along these limits. For the sake of simplicity, the $\omega$ dependency will be explicitly extracted, that is, the limits will be expressed in terms of $\omega \rightarrow 0$.

By following the same reasoning as above, after evaluating the low-frequency asymptotic behavior of the interior and exterior wavenumber, resulting in the operator scalings in real and imaginary parts, we can estimate the asymptotic scalings of the quasi-Helmholtz decomposed blocks of the PMCHWT equation by enforcing the proper cancellations (equations \eqref{eqn:Tpcanc}, \eqref{eqn:K0canc}). 

In addition, notice that the imaginary part of $\mat T_{A,k}$ may present a different scaling depending on whether the source or test basis function is solenoidal or non-solenoidal. Indeed, in the case of solenoidal basis functions, the second term in the Green's function Maclaurin expansion vanishes.

Then, from the above mentioned decomposition, the asymptotic scalings of the quasi-Helmholtz decomposed inverse matrix $\mat Z_{\mat \Lambda \mat H \mat \Sigma}^{-1}$ can be retrieved, for example by applying the Woodbury formula \cite{henderson1981deriving}. These results are summarized in the following.

\subsection{The Quasi-Static Regime}
As already shown in \Cref{sec:analysisQSR}, the interior wavenumber asymptotic behavior in the quasi-static regime is $\Re\{k_1 L\}=\mathcal{O}(\omega)$, $\Im\{k_1 L\}=\mathcal{O}(\omega)$.
This results in the loop-star decomposition of matrix $\bar{\mat Z}$ in the low-frequency limit as
\begin{align}
    \Re\{\mat Z_{\mat \Lambda \mat H \mat \Sigma}\}&=\mathcal{O} 
    \begin{pmatrix}
 \omega^3 & \omega^3 & \omega^3 &  
 \omega^2 &  \omega^2 & 1 \\
 \omega^3 & \omega^3 & \omega^3 &  
 \omega^2 &  1 & 1\\
 \omega^3 & \omega^3 & \omega^{-1} &  
 1 &  1 & 1  \\
 \omega^2 & \omega^2 & 1 & 
 \omega & \omega &  \omega \\
 \omega^2 & 1 & 1 & 
 \omega & \omega &  \omega \\
 1 & 1 & 1 & 
 \omega & \omega &  1 
    \end{pmatrix}\\
    \Im\{\mat Z_{\mat \Lambda \mat H \mat \Sigma}\}&=\mathcal{O} 
    \begin{pmatrix}
 \omega & \omega & \omega &  
 \omega^2 &  \omega^2 & \omega^2 \\
 \omega & \omega & \omega &  
 \omega^2 &  \omega^2 & \omega^2\\
 \omega & \omega & \omega^{-1} &  
 \omega^2 &  \omega^2 & \omega^2  \\
 \omega^2 & \omega^2 & \omega^2 & 
 \omega & \omega &  \omega \\
 \omega^2 & \omega^2 & \omega^2 & 
 \omega & \omega &  \omega \\
 \omega^2 & \omega^2 & \omega^2 & 
 \omega & \omega &  \omega^{-1} 
    \end{pmatrix}\,.
\end{align}
The loop-star decomposition of the matrix inverse behaves as
\begin{align}
    \Re\{\mat Z_{\mat \Lambda \mat H \mat \Sigma}^{-1}\}&=\mathcal{O} 
    \begin{pmatrix}
 1 & \omega & \omega &  
 1 &  1 & 1 \\
 \omega & \omega & \omega &  
 1 &  1 & \omega^2\\
 \omega & \omega & \omega &  
 1 &  \omega^2 & \omega^2  \\
 1 & 1 & 1 & 
 \omega^{-1} & \omega &  \omega \\
 1 & 1 & \omega^2 & 
 \omega & \omega^2 &  \omega^2 \\
 1 & \omega^2 & \omega^2 & 
 \omega & \omega^2 &  \omega^2 
    \end{pmatrix}\\
    \Im\{\mat Z_{\mat \Lambda \mat H \mat \Sigma}^{-1}\}&=\mathcal{O} 
    \begin{pmatrix}
 \omega^{-1} & \omega & \omega &  
 1 &  \omega & \omega \\
 \omega & \omega & \omega &  
 1 &  \omega^2 & \omega^2\\
 \omega & \omega & \omega &  
 1 &  \omega^2 & \omega^2  \\
 1 & 1 & 1 & 
 \omega^{-1} & \omega &  \omega \\
 \omega & \omega^2 & \omega^2 & 
 \omega & \omega &  \omega \\
 \omega & \omega^2 & \omega^2 & 
 \omega & \omega &  \omega 
    \end{pmatrix}\,.
\end{align}

\subsection{The Eddy-Current-Free Eddy-Current Regime}
The low-frequency asymptotic behavior of the interior wavenumber in the ECFR is $\Re\{k_1 L\}=\mathcal{O}(\omega^{1/2})$, $\Im\{k_1 L\}=\mathcal{O}(\omega^{1/2})$ for $\omega \rightarrow 0$.
This results in the following asymptotic behavior of the loop-star decomposed system and its inverse
\begin{align}
    \Re\{\mat Z_{\mat \Lambda \mat H \mat \Sigma}\}&=\mathcal{O} 
    \begin{pmatrix}
 \omega^2 & \omega^2 & \omega^2 &  
 \omega^{3/2} &  \omega^{3/2} & 1 \\
 \omega^2 & \omega^2 & \omega^2 &  
 \omega^{3/2} &  1 & 1\\
 \omega^2 & \omega^2 & 1 &  
 1 &  1 & 1  \\
 \omega^{3/2} & \omega^{3/2} & 1 & 
 1 & 1 &  1 \\
 \omega^{3/2} & 1 & 1 & 
 1 & 1 &  1 \\
 1 & 1 & 1 & 
 1 & 1 &  \omega^{-1/2}
    \end{pmatrix}\\
    \Im\{\mat Z_{\mat \Lambda \mat H \mat \Sigma}\}&=\mathcal{O} 
    \begin{pmatrix}
 \omega & \omega & \omega &  
 \omega &  \omega & \omega \\
 \omega & \omega & \omega &  
 \omega &  \omega & \omega\\
 \omega & \omega & \omega^{-1} &  
 \omega &  \omega & \omega  \\
 \omega & \omega & \omega & 
 \omega & \omega &  \omega \\
 \omega & \omega & \omega & 
 \omega & \omega &  \omega \\
 \omega & \omega & \omega & 
 \omega & \omega &  \omega^{-1} 
    \end{pmatrix}
\end{align}
and
\begin{align}
    \Re\{\mat Z_{\mat \Lambda \mat H \mat \Sigma}^{-1}\}&=\mathcal{O} 
    \begin{pmatrix}
 \omega^{-1/2} & 1 & \omega^{3/2} &  
 1 &  1 & 1 \\
 1 & 1 & \omega^2 &  
 1 &  1 & \omega^{3/2}\\
 \omega^{3/2} & \omega^2 & \omega^2 &  
 \omega^2 &  \omega^2 & \omega^2  \\
 1 & 1 & \omega^2 & 
 1 & \omega^{3/2} &  \omega^{3/2} \\
 1 & 1 & \omega^2 & 
 \omega^{3/2} & \omega^{3/2} &  \omega^{3/2} \\
 1 & \omega^{3/2} & \omega^2 & 
 \omega^{3/2} & \omega^{3/2} &  \omega^{3/2} 
    \end{pmatrix}\\
    \Im\{\mat Z_{\mat \Lambda \mat H \mat \Sigma}^{-1}\}&=\mathcal{O} 
    \begin{pmatrix}
 \omega^{-1} & \omega^{1/2} & \omega &  
 \omega^{1/2} &  \omega^{1/2} & \omega^{1/2} \\
 \omega^{1/2} & \omega & \omega &  
 \omega &  \omega & \omega\\
 \omega & \omega & \omega &  
 \omega &  \omega^{5/2} & \omega^{5/2}  \\
 \omega^{1/2} & \omega & \omega & 
 \omega & \omega &  \omega \\
 \omega^{1/2} & \omega & \omega^{5/2} & 
 \omega & \omega &  \omega \\
 \omega^{1/2} & \omega & \omega^{5/2} & 
 \omega & \omega &  \omega 
    \end{pmatrix}
\end{align}

\subsection{The Skin-Effect-Dominated Eddy-Current Regime}
The product $(k_1L)$ stays constant in the low-frequency limit considered in this case, that is, $\Re\{k_1 L\}=\mathcal{O}(1)$, $\Im\{k_1 L\}=\mathcal{O}(1)$ for $\omega \rightarrow 0$.
The loop-star decomposition of matrix $\bar{\mat Z}$ then behaves as
\begin{align}
    \Re\{\mat Z_{\mat \Lambda \mat H \mat \Sigma}\}&=\mathcal{O} 
    \begin{pmatrix}
 \omega & \omega & \omega &  
 1 &  1 & 1 \\
 \omega & \omega & \omega &  
 1 &  1 & 1\\
 \omega & \omega & 1 &  
 1 &  1 & 1  \\
 1 & 1 & 1 & 
 \omega^{-1} & \omega^{-1} &  \omega^{-1} \\
 1 & 1 & 1 & 
 \omega^{-1} & \omega^{-1} &  \omega^{-1} \\
 1 & 1 & 1 & 
 \omega^{-1} & \omega^{-1} &  \omega^{-1} 
    \end{pmatrix}\\
    \Im\{\mat Z_{\mat \Lambda \mat H \mat \Sigma}\}&=\mathcal{O} 
    \begin{pmatrix}
 \omega & \omega & \omega &  
 1 &  1 & 1 \\
 \omega & \omega & \omega &  
 1 &  1 & 1\\
 \omega & \omega & \omega^{-1} &  
 1 &  1 & 1  \\
 1 & 1 & 1 & 
 \omega^{-1} & \omega^{-1} &  \omega^{-1} \\
 1 & 1 & 1 & 
 \omega^{-1} & \omega^{-1} &  \omega^{-1} \\
 1 & 1 & 1 & 
 \omega^{-1} & \omega^{-1} &  \omega^{-1} 
    \end{pmatrix}
\end{align}
The loop-star decomposition of the matrix inverse is
\begin{align}
    \Re\{\mat Z^{-1}_{\mat \Lambda \mat H \mat \Sigma}\}&=\mathcal{O} 
    \begin{pmatrix}
 \omega^{-1} & \omega^{-1} & \omega &  
 1 &  1 & 1 \\
 \omega^{-1} & \omega^{-1} & \omega &  
 1 &  1 & 1\\
 \omega & \omega & \omega^2 &  
 \omega^2 & \omega^2 & \omega^2  \\
 1 & 1 & \omega^2 & 
 \omega & \omega &  \omega \\
 1 & 1 & \omega^2 & 
 \omega & \omega &  \omega \\
 1 & 1 & \omega^2 & 
 \omega & \omega &  \omega 
    \end{pmatrix}\\
    \Im\{\mat Z^{-1}_{\mat \Lambda \mat H \mat \Sigma}\}&=\mathcal{O} 
    \begin{pmatrix}
 \omega^{-1} & \omega^{-1} & \omega &  
 1 &  1 & 1 \\
 \omega^{-1} & \omega^{-1} & \omega &  
 1 &  1 & 1\\
 \omega & \omega & \omega & 
 \omega^2 & \omega^2 & \omega^2 \\
 1 & 1 & \omega^2 & 
 \omega & \omega &  \omega \\
 1 & 1 & \omega^2 & 
 \omega & \omega &  \omega \\
 1 & 1 & \omega^2 & 
 \omega & \omega &  \omega 
    \end{pmatrix}
\end{align}

\subsection{Excitation Definition and Resulting Solution}
\label{sec:dominantexcitation}
The behavior of the different components of the excitation vector are required to determine the behavior of the components of the solution of the PMCHWT. In this work, we focus on the magnetic frill excitation \cite{tsai1972numerical,popovic1982analysis}, as one of the most common ways to model an imposed potential difference in the circuital framework. Notice, however, that the same procedure used in this work to define a preconditioning strategy for the PMCHWT equation optimized for  magnetic frill sources can also be applied to other types of excitations.

The low-frequency asymptotic scalings of the quasi-Helmholtz components of the electric and magnetic fields produced by a magnetic frill have already been studied in \cite{chhim2020eddy} and are reported in the following tables. It is worth noticing that, following the reasoning in \cite{bogaert2011low}, the harmonic component of the electric field shows a different behavior depending on whether there exists a global loop of the geometry passing through the magnetic frill or not, resulting in an inductive excitation or a capacitive one respectively, both of which will be studied in the following.

By applying the inverse of the operator matrices to the resulting RHSs, it is possible to retrieve the low-frequency asymptotic scalings of the quasi-Helmholtz components of the surface current densities. The scattering of each of them, in their real and imaginary parts, results in a component of the electromagnetic fields.
As a result of this analysis, the current components required to correctly reconstruct the fields in the three regimes for the two types of excitations considered are summarized in tables \ref{table:indtablecurrent} and \ref{table:captablecurrent}, in the right column of sections (c), (d), (e), and (f).

When solving numerically the standard PMCHWT in floating point arithmetic, the real parts (or imaginary parts) of the loop and star components of the currents are stored in the same floating point number. In the low-frequency limit, this determines a loss of information about the non-dominant current components, lying outside the dynamic range of the floating point number determined by the larger current contributions. Nevertheless, the lost components may be crucial in the reconstruction of the fields. Indeed, by comparing the current components required for the correct evaluation of the scattered fields with the ones that can be retrieved without proper rescaling of the quasi-Helmholtz components of the system in the right column of section (b), one can notice that in many cases the current resulting from a naive solution of the PMCHWT equation does not allow for the correct reconstruction of the fields. For example, from Table~\ref{table:indtablecurrent} we infer that, in the EDCR, the current component $\Im(\veg j_{\mSig})$ is required for a correct reconstruction of the exterior electric field; however this component is lost in numerical cancellation in the low-frequency limit if a proper rescaling of the current coefficients is not applied, resulting in a deterioration of the accuracy of the scattered field.

Now that the critical components of the current have been identified, we can design our preconditioning strategy with the objective of preserving them in addition to curing the ill-conditioning of the matrix.
\begin{table*}[!t]
    \renewcommand{\arraystretch}{1.4}
    \setlength\tabcolsep{5pt}
    \caption{Scalings of the real and imaginary parts of quantities of interest when $\omega \to 0$ for an inductive magnetic frill excitation }
    \label{table:indtablecurrent}
    \centering
    \begin{tabular}{|l|ccc|ccc|c|}
        \hline
        \multicolumn{8}{|c|}{\textbf{(a) Right hand side}} \\ 
        \hline
         Excitation & $(\Re,\Im)(\veg E^i_{\mLam})$ & $(\Re,\Im)(\veg E^i_{\mat H})$ & $(\Re,\Im)(\veg E^i_{\mSig})$ & $(\Re,\Im)(\veg H^i_{\mLam})$ & $(\Re,\Im)(\veg H^i_{\mat H})$ & $(\Re,\Im)(\veg H^i_{\mSig})$ &  \\
        \hline
        {Inductive} & ($\omega^2$, $\omega^3$) & ($1$, $\omega^3$) & (1, $\omega^3$) & ($\omega^4$, $\omega$) & ($\omega^4$, $\omega$) & ($\omega^4$, $\omega$) &  \\
        \hline
        \multicolumn{8}{|c|}{\textbf{(b) Surface current density}} \\
        \hline
        Regime & $(\Re,\Im)(\veg j_{\mLam})$ & $(\Re,\Im)(\veg j_{\mat H})$ & $(\Re,\Im)(\veg j_{\mSig})$ & $(\Re,\Im)(\veg m_{\mLam})$ & $(\Re,\Im)(\veg m_{\mat H})$ & $(\Re,\Im)(\veg m_{\mSig})$ & Current dominant components \\
        \hline
        {QSR} & (\blue{$\omega$}, \blue{$\omega$}) & (\blue{$\omega$}, $\omega$) & ($\omega$, $\omega$) & (\blue{1}, \blue{1}) & (\blue{1}, \blue{$\omega^2$}) & ($\omega^{2}$, \blue{$\omega^2$}) & \makecell{$\Re(\veg j_{\mLam},\veg j_{\mat H},\veg j_{\mSig},\veg m_{\mLam},\veg m_{\mat H})$\\$\Im(\veg j_{\mLam},\veg j_{\mat H},\veg j_{\mSig},\veg m_{\mLam})$} \\
        \hline
        {ECFR} & (\blue{$1$}, \blue{$\omega^{1/2}$}) & (\blue{$1$}, $\omega$) & ($\omega^2$, $\omega$) & (\blue{1}, $\omega$) & (\blue{1}, \blue{$\omega$}) & ($\omega^{3/2}$, \blue{$\omega$}) & \makecell{$\Re(\veg j_{\mLam},\veg j_{\mat H},\veg m_{\mLam},\veg m_{\mat H})$\\$\Im(\veg j_{\mLam},\veg m_{\mLam},\veg m_{\mat H},\veg m_{\mSig})$} \\
        \hline
        {SEDR} & (\blue{$\omega^{-1}$}, \blue{$\omega^{-1}$}) & (\blue{$\omega^{-1}$}, $\omega^{-1}$) & ($\omega$, $\omega$) & (\blue{1}, $1$) & (\blue{1}, \blue{$1$}) & ($1$, \blue{$1$}) & \makecell{$\Re(\veg j_{\mLam},\veg j_{\mat H},\veg m_{\mLam},\veg m_{\mat H},\veg m_{\mSig})$\\$\Im(\veg j_{\mLam},\veg j_{\mat H},\veg m_{\mLam},\veg m_{\mat H},\veg m_{\mSig})$} \\
        \hline
        
        \multicolumn{8}{|c|}{\textbf{(c) Electric interior near field}} \\ 
        \hline
        Regime & \makecell{$\veg E(\Re \veg j_{\mLam})$\\$\veg E(\Im \veg j_{\mLam})$} & \makecell{$\veg E(\Re \veg j_{\mat H})$\\$\veg E(\Im \veg j_{\mat H})$} & \makecell{$\veg E(\Re \veg j_{\mSig})$\\$\veg E(\Im \veg j_{\mSig})$} & \makecell{$\veg E(\Re \veg m_{\mLam})$\\$\veg E(\Im \veg m_{\mLam})$} & \makecell{$\veg E(\Re \veg m_{\mat H})$\\$\veg E(\Im \veg m_{\mat H})$} & \makecell{$\veg E(\Re \veg m_{\mSig})$\\$\veg E(\Im \veg m_{\mSig})$} & Current dominant components \\
        \hline
        QSR &\makecell{ ({$\omega^4$}, $\omega^2$)\\({$\omega^{2}$}, \blue{$\omega^{4}$}) }&\makecell{ ($\omega^{4}$, $\omega^{2}$)\\ ($\omega^2$, $\omega^4$)}&\makecell{ ($1$, $1$)\\({$1$}, $1$) }&\makecell{ (\blue{$1$}, $\omega^{2}$)\\($\omega^2$, $1$) }&\makecell{ (\blue{$1$}, $\omega^{2}$)\\($\omega^4$, $\omega^2$) }&\makecell{ ($\omega^{2}$, {$\omega^{4}$})\\({$\omega^4$}, $\omega^2$) }&\makecell{$\Re( \veg j_{\mSig},  \veg m_{\mLam},  \veg m_{\mat H})$\\$\Im(\veg j_{\mSig},\veg m_{\mLam})$} \\
        \hline
        ECFR &\makecell{ ({$\omega^2$}, $\omega$)\\({$\omega^{3/2}$}, \blue{$\omega^{5/2}$}) }&\makecell{ ($\omega^{2}$, $\omega$)\\ ($\omega^2$, $\omega^3$)}&\makecell{ ($\omega^{2}$, $\omega^3$)\\ ($\omega^2$, $\omega$)}&\makecell{ (\blue{$1$}, $\omega$)\\($\omega^2$, $\omega$) }&\makecell{ (\blue{$1$}, $\omega$)\\($\omega^2$, $\omega$) }&\makecell{ ($\omega^{3/2}$, {$\omega^{5/2}$})\\({$\omega^2$}, $\omega$) }&$\Re(   \veg m_{\mLam},  \veg m_{\mat H})$ \\
        \hline
        SEDR &\makecell{ ({$1$}, $1$)\\({$1$}, \blue{$1$}) }&\makecell{ ($1$, $1$)\\ ($1$, $1$)}&\makecell{ ($\omega^2$, $\omega^2$)\\({$\omega^2$}, $\omega^2$) }&\makecell{ ({$1$}, $1$)\\({$1$}, \blue{$1$}) }&\makecell{ ({$1$}, $1$)\\({$1$}, \blue{$1$}) }&\makecell{ ({$1$}, $1$)\\({$1$}, \blue{$1$}) }&\makecell{$\Re( \veg j_{\mLam},\veg j_{\mat H},  \veg m_{\mLam},  \veg m_{\mat H},\veg m_{\mSig})$\\$\Im( \veg j_{\mLam},\veg j_{\mat H},  \veg m_{\mLam},  \veg m_{\mat H},\veg m_{\mSig})$} \\
        \hline

        \multicolumn{8}{|c|}{\textbf{(d) Magnetic interior near field}} \\ 
        \hline
        Regime & \makecell{$\veg H(\Re \veg j_{\mLam})$\\$\veg H(\Im \veg j_{\mLam})$} & \makecell{$\veg H(\Re \veg j_{\mat H})$\\$\veg H(\Im \veg j_{\mat H})$} & \makecell{$\veg H(\Re \veg j_{\mSig})$\\$\veg H(\Im \veg j_{\mSig})$} & \makecell{$\veg H(\Re \veg m_{\mLam})$\\$\veg H(\Im \veg m_{\mLam})$} & \makecell{$\veg H(\Re \veg m_{\mat H})$\\$\veg H(\Im \veg m_{\mat H})$} & \makecell{$\veg H(\Re \veg m_{\mSig})$\\$\veg H(\Im \veg m_{\mSig})$} & Current dominant components \\
        \hline
        QSR  &\makecell{ ({$\omega$}, $\omega^3$)\\({$\omega^{3}$}, \blue{$\omega$}) }&\makecell{ ($\omega$, $\omega^{3}$)\\ ($\omega^3$, $\omega$)}&\makecell{ ($\omega$, $\omega^3$)\\({$\omega^3$}, $\omega$) }&\makecell{ (\blue{$\omega$}, $\omega$)\\($\omega$, $\omega$) }&\makecell{ (\blue{$\omega$}, $\omega$)\\($\omega^3$, $\omega^3$) }&\makecell{ ($\omega^{3}$, {$\omega$})\\({$\omega$}, $\omega^3$) }&\makecell{$\Re({\veg j_{\mLam}},{\veg j_{\mat H}}, \veg j_{\mSig},  \veg m_{\mLam},  \veg m_{\mat H},{\veg m_{\mSig}})$\\$\Im({\veg j_{\mLam}},{\veg j_{\mat H}},\veg j_{\mSig},\veg m_{\mLam},{\veg m_{\mSig}})$}  \\
        
        \hline
        ECFR &\makecell{ ({$1$}, $\omega$)\\({$\omega^{3/2}$}, \blue{$\omega^{1/2}$}) }&\makecell{ ($1$, $\omega$)\\ ($\omega^2$, $\omega$)}&\makecell{ ($\omega^2$, $\omega^3$)\\({$\omega^2$}, $\omega$) }&\makecell{ (\blue{$1$}, $\omega$)\\($\omega^2$, $\omega$) }&\makecell{ (\blue{$1$}, $\omega$)\\($\omega^2$, $\omega$) }&\makecell{ ($\omega^{3/2}$, {$\omega^{1/2}$})\\({$1$}, $\omega$) }&\makecell{$\Re(\veg j_{\mLam},\veg j_{\mat H},  \veg m_{\mLam},  \veg m_{\mat H})$\\$\Im(\veg m_{\mSig})$}  \\
        \hline
        SEDR &\makecell{ ({$\omega^{-1}$}, $\omega^{-1}$)\\({$\omega^{-1}$}, \blue{$\omega^{-1}$}) }&\makecell{ ($\omega^{-1}$, $\omega^{-1}$)\\ ($\omega^{-1}$, $\omega^{-1}$)}&\makecell{ ($\omega$, $\omega$)\\({$\omega$}, $\omega$) }&\makecell{ (\blue{$\omega^{-1}$}, $\omega^{-1}$)\\($\omega^{-1}$, $\omega^{-1}$) }&\makecell{ (\blue{$\omega^{-1}$}, $\omega^{-1}$)\\($\omega^{-1}$, $\omega^{-1}$) }&\makecell{ ($\omega^{-1}$, {$\omega^{-1}$})\\({$\omega^{-1}$}, $\omega^{-1}$) }&\makecell{$\Re(\veg j_{\mLam},\veg j_{\mat H},  \veg m_{\mLam},  \veg m_{\mat H},\veg m_{\mSig})$\\$\Im(\veg j_{\mLam},\veg j_{\mat H},  \veg m_{\mLam},  \veg m_{\mat H},\veg m_{\mSig})$}  \\

        \hline

        \multicolumn{8}{|c|}{\textbf{(e) Electric exterior near field}} \\ 
        \hline
        Regime & \makecell{$\veg E(\Re \veg j_{\mLam})$\\$\veg E(\Im \veg j_{\mLam})$} & \makecell{$\veg E(\Re \veg j_{\mat H})$\\$\veg E(\Im \veg j_{\mat H})$} & \makecell{$\veg E(\Re \veg j_{\mSig})$\\$\veg E(\Im \veg j_{\mSig})$} & \makecell{$\veg E(\Re \veg m_{\mLam})$\\$\veg E(\Im \veg m_{\mLam})$} & \makecell{$\veg E(\Re \veg m_{\mat H})$\\$\veg E(\Im \veg m_{\mat H})$} & \makecell{$\veg E(\Re \veg m_{\mSig})$\\$\veg E(\Im \veg m_{\mSig})$} & Current dominant components \\
        \hline
        QSR  &\makecell{ ({$\omega^5$}, $\omega^2$)\\({$\omega^{2}$}, \blue{$\omega^5$}) }&\makecell{ ($\omega^5$, $\omega^{2}$)\\ ($\omega^2$, $\omega^5$)}&\makecell{ ($\omega^3$, $1$)\\({$1$}, $\omega^3$) }&\makecell{ (\blue{$1$}, $\omega^3$)\\($\omega^3$, $1$) }&\makecell{ (\blue{$1$}, $\omega^3$)\\($\omega^5$, $\omega^2$) }&\makecell{ ($\omega^{2}$, {$\omega^5$})\\({$\omega^5$}, $\omega^2$) }&\makecell{$\Re( \veg j_{\mSig},  \veg m_{\mLam},  \veg m_{\mat H})$\\$\Im(\veg j_{\mSig},\veg m_{\mLam})$}  \\
        \hline
        ECFR  &\makecell{ ({$\omega^4$}, $\omega$)\\({$\omega^{3/2}$}, \blue{$\omega^{9/2}$}) }&\makecell{ ($\omega^4$, $\omega$)\\ ($\omega^2$, $\omega^5$)}&\makecell{ ($\omega^4$, $\omega$)\\({$1$}, $\omega^3$) }&\makecell{ (\blue{$1$}, $\omega^3$)\\($\omega^4$, $\omega$) }&\makecell{ (\blue{$1$}, $\omega^3$)\\($\omega^4$, $\omega$) }&\makecell{ ($\omega^{3/2}$, {$\omega^{9/2}$})\\({$\omega^4$}, $\omega$) }&\makecell{$\Re( \veg m_{\mLam},  \veg m_{\mat H})$\\$\Im(\veg j_{\mSig})$}  \\
        \hline
        SEDR  &\makecell{ ({$\omega^3$}, $1$)\\({$1$}, \blue{$\omega^{3}$}) }&\makecell{ ({$\omega^3$}, $1$)\\({$1$}, \blue{$\omega^{3}$})}&\makecell{ ({$\omega^3$}, $1$)\\({$1$}, \blue{$\omega^{3}$}) }&\makecell{ (\blue{$1$}, $\omega^3$)\\($\omega^3$, $1$) }&\makecell{ (\blue{$1$}, $\omega^3$)\\($\omega^3$, $1$) }&\makecell{ (\blue{$1$}, $\omega^3$)\\($\omega^3$, $1$) }&\makecell{$\Re( \veg j_{\mLam},\veg j_{\mat H},\veg j_{\mSig},\veg m_{\mLam},\veg m_{\mat H}, \veg m_{\mSig})$ \\$\Im( \veg j_{\mLam},\veg j_{\mat H},\veg j_{\mSig},\veg m_{\mLam},\veg m_{\mat H}, \veg m_{\mSig})$}  \\
        \hline

        \multicolumn{8}{|c|}{\textbf{(f) Magnetic exterior near field}} \\ 
        \hline
        Regime & \makecell{$\veg H(\Re \veg j_{\mLam})$\\$\veg H(\Im \veg j_{\mLam})$} & \makecell{$\veg H(\Re \veg j_{\mat H})$\\$\veg H(\Im \veg j_{\mat H})$} & \makecell{$\veg H(\Re \veg j_{\mSig})$\\$\veg H(\Im \veg j_{\mSig})$} & \makecell{$\veg H(\Re \veg m_{\mLam})$\\$\veg H(\Im \veg m_{\mLam})$} & \makecell{$\veg H(\Re \veg m_{\mat H})$\\$\veg H(\Im \veg m_{\mat H})$} & \makecell{$\veg H(\Re \veg m_{\mSig})$\\$\veg H(\Im \veg m_{\mSig})$} & Current dominant components \\
        \hline
        QSR  &\makecell{ ({$\omega$}, $\omega^4$)\\({$\omega^{4}$}, \blue{$\omega$}) }&\makecell{ ($\omega$, $\omega^{4}$)\\ ($\omega^4$, $\omega$)}&\makecell{ ($\omega$, $\omega^4$)\\({$\omega^4$}, $\omega$) }&\makecell{ (\blue{$\omega^4$}, $\omega$)\\($\omega$, $\omega^4$) }&\makecell{ (\blue{$\omega^4$}, $\omega$)\\($\omega^3$, $\omega^6$) }&\makecell{ ($\omega^{4}$, {$\omega$})\\({$\omega$}, $\omega^4$) }&\makecell{$\Re( {\veg j_{\mLam}},{\veg j_{\mat H}},\veg j_{\mSig},\veg m_{\mLam},\veg m_{\mat H}, {\veg m_{\mSig}})$ \\$\Im( {\veg j_{\mLam}},{\veg j_{\mat H}},\veg j_{\mSig},\veg m_{\mLam}, {\veg m_{\mSig}})$}  \\
        \hline
        ECFR  &\makecell{ ({$1$}, $\omega^3$)\\({$\omega^{7/2}$}, \blue{$\omega^{1/2}$}) }&\makecell{ ($1$, $\omega^{3}$)\\ ($\omega^4$, $\omega$)}&\makecell{ ($\omega^2$, $\omega^5$)\\({$\omega^4$}, $\omega$) }&\makecell{ (\blue{$\omega^4$}, $\omega$)\\($\omega^2$, $\omega^5$) }&\makecell{ (\blue{$\omega^4$}, $\omega$)\\($\omega^2$, $\omega^5$) }&\makecell{ ($\omega^{7/2}$, {$\omega^{1/2}$})\\({$1$}, $\omega^3$) }&\makecell{$\Re( {\veg j_{\mLam}},{\veg j_{\mat H}})$ \\$\Im(  {\veg m_{\mSig}})$}  \\
        \hline
        SEDR  &\makecell{ ({$\omega^{-1}$}, $\omega^2$)\\({$\omega^{2}$}, \blue{$\omega^{-1}$}) }&\makecell{ ({$\omega^{-1}$}, $\omega^2$)\\({$\omega^{2}$}, \blue{$\omega^{-1}$}) }&\makecell{ ($\omega$, $\omega^4$)\\({$\omega^4$}, $\omega$) }&\makecell{ (\blue{$\omega^4$}, $\omega$)\\($\omega$, $\omega^4$) }&\makecell{ (\blue{$\omega^4$}, $\omega$)\\($\omega$, $\omega^4$) }&\makecell{ ($\omega^{2}$, {$\omega^{-1}$})\\({$\omega^{-1}$}, $\omega^2$) }&\makecell{$\Re( {\veg j_{\mLam}},{\veg j_{\mat H}},\veg m_{\mSig})$ \\$\Im( {\veg j_{\mLam}},{\veg j_{\mat H}},\veg m_{\mSig})$}  \\
        \hline

        \multicolumn{8}{|c|}{\textbf{(g) Rescaled current density}} \\ 
        \hline
        Regime & $a_R^{-1} \veg j_{\mLam}$ & $a_R^{-1} \veg j_{\mat H}$ & $b_R^{-1} \veg j_{\mSig}$ & $c_R^{-1} \veg m_{\mLam}$ & $c_R^{-1} \veg m_{\mat H}$ & $d_R^{-1} \veg m_{\mSig}$ & Recovered components \\
        \hline
        QSR & ($1,1$) & ($1,1$) & ($1,1$) & ($1,1$) & ($1, \omega^{2}$) & ($1,1$) & \makecell{$\Re( \veg j_{\mLam},\veg j_{\mat H},\veg j_{\mSig}, \veg m_{\mLam},\veg m_{\mat H},\veg m_{\mSig})$ \\$\Im( \veg j_{\mLam},\veg j_{\mat H},\veg j_{\mSig}, \veg m_{\mLam},\veg m_{\mSig})$} \\
        \hline
        ECFR & ($\omega^{1/2}, \omega$) & ($\omega^{1/2}, \omega^{3/2}$) & ($\omega, 1$) & ($1, \omega$) & ($1, \omega$) & ($\omega, \omega^{1/2}$) & \makecell{$\Re( \veg j_{\mLam},\veg j_{\mat H}, \veg m_{\mLam},\veg m_{\mat H})$ \\$\Im(\veg j_{\mSig}, \veg m_{\mSig})$} \\
        \hline
        SEDR & ($\omega^{-1/2}, \omega^{-1/2}$) & ($\omega^{-1/2}, \omega^{-1/2}$) & ($\omega^{-1/2}, \omega^{-1/2}$) & ($\omega^{-1/2}, \omega^{-1/2}$) & ($\omega^{-1/2}, \omega^{-1/2}$) & ($\omega^{-1/2}, \omega^{-1/2}$) & \makecell{$\Re( \veg j_{\mLam},\veg j_{\mat H},\veg j_{\mSig},\veg m_{\mLam},\veg m_{\mat H}, \veg m_{\mSig})$ \\$\Im( \veg j_{\mLam},\veg j_{\mat H},\veg j_{\mSig},\veg m_{\mLam},\veg m_{\mat H}, \veg m_{\mSig})$} \\
        \hline
    \end{tabular}
\end{table*}
\begin{table*}[!t]
    \renewcommand{\arraystretch}{1.4}
    \setlength\tabcolsep{5pt}
    \caption{Scalings of the real and imaginary parts of quantities of interest when $\omega \to 0$ for a capacitive magnetic frill excitation}
    \label{table:captablecurrent}
    \centering
    \begin{tabular}{|l|ccc|ccc|c|}
        \hline
        \multicolumn{8}{|c|}{\textbf{(a) Right hand side}} \\ 
        \hline
         Excitation & $(\Re,\Im)(\veg E^i_{\mLam})$ & $(\Re,\Im)(\veg E^i_{\mat H})$ & $(\Re,\Im)(\veg E^i_{\mSig})$ & $(\Re,\Im)(\veg H^i_{\mLam})$ & $(\Re,\Im)(\veg H^i_{\mat H})$ & $(\Re,\Im)(\veg H^i_{\mSig})$ &  \\
        \hline
        {Capacitive} & ($\omega^2$, $\omega^3$) & ($\omega^2$, $\omega^3$) & (1, $\omega^3$) & ($\omega^4$, $\omega$) & ($\omega^4$, $\omega$) & ($\omega^4$, $\omega$) &  \\
        \hline
        \multicolumn{8}{|c|}{\textbf{(b) Surface current density}} \\
        \hline
        Regime & $(\Re,\Im)(\veg j_{\mLam})$ & $(\Re,\Im)(\veg j_{\mat H})$ & $(\Re,\Im)(\veg j_{\mSig})$ & $(\Re,\Im)(\veg m_{\mLam})$ & $(\Re,\Im)(\veg m_{\mat H})$ & $(\Re,\Im)(\veg m_{\mSig})$ & Current dominant components \\
        \hline
        {QSR} & (\blue{$\omega$}, \blue{$\omega$}) & (\blue{$\omega$}, $\omega$) & ($\omega$, $\omega$) & (\blue{1}, \blue{1}) & (\blue{$\omega^2$}, \blue{$\omega^2$}) & ($\omega^{2}$, \blue{$\omega^2$}) & \makecell{$\Re(\veg j_{\mLam},\veg j_{\mat H},\veg j_{\mSig},\veg m_{\mLam})$\\$\Im(\veg j_{\mLam},\veg j_{\mat H},\veg j_{\mSig},\veg m_{\mLam})$} \\
        \hline
        {ECFR} & (\blue{$\omega^{3/2}$}, \blue{$\omega$}) & (\blue{$\omega^{2}$}, $\omega$) & ($\omega^2$, $\omega$) & (\blue{$\omega^2$}, $\omega$) & (\blue{$\omega^2$}, \blue{$\omega^{5/2}$}) & ($\omega^{2}$, \blue{$\omega^{5/2}$}) & \makecell{$\Re(\veg j_{\mLam},\veg m_{\mLam},\veg m_{\mat H},\veg m_{\mSig})$\\$\Im(\veg j_{\mLam},\veg j_{\mat H},\veg j_{\mSig},\veg m_{\mLam})$} \\
        \hline
        {SEDR} & (\blue{$\omega$}, \blue{$\omega$}) & (\blue{$\omega$}, $\omega$) & ($\omega^2$, $\omega$) & (\blue{$\omega^2$}, $\omega^2$) & (\blue{$\omega^2$}, $\omega^2$) & (\blue{$\omega^2$}, $\omega^2$) & \makecell{$\Re(\veg j_{\mLam},\veg j_{\mat H},\veg m_{\mLam},\veg m_{\mat H},\veg m_{\mSig})$\\$\Im(\veg j_{\mLam},\veg j_{\mat H},\veg j_{\mSig},\veg m_{\mLam},\veg m_{\mat H},\veg m_{\mSig})$} \\

        \hline
        \multicolumn{8}{|c|}{\textbf{(c) Electric interior near field}} \\ 
        \hline
        Regime & \makecell{$\veg E(\Re \veg j_{\mLam})$\\$\veg E(\Im \veg j_{\mLam})$} & \makecell{$\veg E(\Re \veg j_{\mat H})$\\$\veg E(\Im \veg j_{\mat H})$} & \makecell{$\veg E(\Re \veg j_{\mSig})$\\$\veg E(\Im \veg j_{\mSig})$} & \makecell{$\veg E(\Re \veg m_{\mLam})$\\$\veg E(\Im \veg m_{\mLam})$} & \makecell{$\veg E(\Re \veg m_{\mat H})$\\$\veg E(\Im \veg m_{\mat H})$} & \makecell{$\veg E(\Re \veg m_{\mSig})$\\$\veg E(\Im \veg m_{\mSig})$} & Current dominant components \\
        \hline
        QSR &\makecell{ ({$\omega^4$}, $\omega^2$)\\({$\omega^{2}$}, \blue{$\omega^{4}$}) }&\makecell{ ($\omega^{4}$, $\omega^{2}$)\\ ($\omega^2$, $\omega^4$)}&\makecell{ ($1$, $1$)\\({$1$}, $1$) }&\makecell{ (\blue{$1$}, $\omega^{2}$)\\($\omega^2$, $1$) }&\makecell{ (\blue{$\omega^{2}$}, $\omega^{4}$)\\($\omega^4$, $\omega^2$) }&\makecell{ ($\omega^{2}$, {$\omega^{4}$})\\({$\omega^4$}, $\omega^2$) }&\makecell{$\Re( \veg j_{\mSig},  \veg m_{\mLam})$\\$\Im(\veg j_{\mSig},\veg m_{\mLam})$} \\
        \hline
        ECFR &\makecell{ ({$\omega^{7/2}$}, $\omega^{5/2}$)\\({$\omega^{2}$}, \blue{$\omega^{3}$}) }&\makecell{ ($\omega^{4}$, $\omega^3$)\\ ($\omega^2$, $\omega^3$)}&\makecell{ ($\omega^{2}$, $\omega^3$)\\ ($\omega^2$, $\omega$)}&\makecell{ (\blue{$\omega^2$}, $\omega^3$)\\($\omega^2$, $\omega$) }&\makecell{ (\blue{$\omega^2$}, $\omega^3$)\\($\omega^{7/2}$, $\omega^{5/2}$) }&\makecell{ ($\omega^{2}$, {$\omega^{3}$})\\({$\omega^{7/2}$}, $\omega^{5/2}$) }&$\Im( \veg j_{\mSig},  \veg m_{\mLam})$ \\
        \hline
        SEDR &\makecell{ ({$\omega^2$}, $\omega^2$)\\({$\omega^2$}, \blue{$\omega^2$}) }&\makecell{ ({$\omega^2$}, $\omega^2$)\\({$\omega^2$}, \blue{$\omega^2$}) }&\makecell{ ($\omega^3$, $\omega^3$)\\({$\omega^2$}, $\omega^2$) }&\makecell{ ({$\omega^2$}, $\omega^2$)\\({$\omega^2$}, \blue{$\omega^2$}) }&\makecell{ ({$\omega^2$}, $\omega^2$)\\({$\omega^2$}, \blue{$\omega^2$}) }&\makecell{ ({$\omega^2$}, $\omega^2$)\\({$\omega^2$}, \blue{$\omega^2$}) }&\makecell{$\Re(\veg j_{\mLam},\veg j_{\mat H}, \veg m_{\mLam},\veg m_{\mat H},\veg m_{\mSig})$ \\$\Im( \veg j_{\mLam},\veg j_{\mat H},\veg j_{\mSig},\veg m_{\mLam},\veg m_{\mat H},\veg m_{\mSig})$} \\
        \hline

        \multicolumn{8}{|c|}{\textbf{(d) Magnetic interior near field}} \\ 
        \hline
        Regime & \makecell{$\veg H(\Re \veg j_{\mLam})$\\$\veg H(\Im \veg j_{\mLam})$} & \makecell{$\veg H(\Re \veg j_{\mat H})$\\$\veg H(\Im \veg j_{\mat H})$} & \makecell{$\veg H(\Re \veg j_{\mSig})$\\$\veg H(\Im \veg j_{\mSig})$} & \makecell{$\veg H(\Re \veg m_{\mLam})$\\$\veg H(\Im \veg m_{\mLam})$} & \makecell{$\veg H(\Re \veg m_{\mat H})$\\$\veg H(\Im \veg m_{\mat H})$} & \makecell{$\veg H(\Re \veg m_{\mSig})$\\$\veg H(\Im \veg m_{\mSig})$} & Current dominant components \\
        \hline
        QSR  &\makecell{ ({$\omega$}, $\omega^3$)\\({$\omega^3$}, \blue{$\omega$}) }&\makecell{ ({$\omega$}, $\omega^3$)\\({$\omega^3$}, \blue{$\omega$}) }&\makecell{ ($\omega$, $\omega^3$)\\({$\omega^3$}, $\omega$) }&\makecell{ ({$\omega$}, $\omega$)\\({$\omega$}, \blue{$\omega$}) }&\makecell{ ({$\omega^3$}, $\omega^3$)\\({$\omega^3$}, \blue{$\omega^3$}) }&\makecell{ ({$\omega^3$}, $\omega$)\\({$\omega$}, \blue{$\omega^3$}) }&\makecell{$\Re({\veg j_{\mLam}},\veg j_{\mat H},\veg j_{\mSig}, \veg m_{\mLam},{\veg m_{\mSig}})$ \\$\Im({\veg j_{\mLam}},\veg j_{\mat H},\veg j_{\mSig}, \veg m_{\mLam},{\veg m_{\mSig}})$}  \\
        \hline
        ECFR  &\makecell{ ({$\omega^{3/2}$}, $\omega^{5/2}$)\\({$\omega^2$}, \blue{$\omega$}) }&\makecell{ ({$\omega^2$}, $\omega^3$)\\({$\omega^2$}, \blue{$\omega$}) }&\makecell{ ($\omega^2$, $\omega^3$)\\({$\omega^2$}, $\omega$) }&\makecell{ ({$\omega^2$}, $\omega^3$)\\({$\omega^2$}, \blue{$\omega$}) }&\makecell{ ({$\omega^2$}, $\omega^3$)\\({$\omega^{7/2}$}, \blue{$\omega^{5/2}$}) }&\makecell{ ({$\omega^2$}, $\omega$)\\({$\omega^{3/2}$}, \blue{$\omega^{5/2}$}) }&\makecell{$\Re(\veg m_{\mSig})$ \\$\Im({\veg j_{\mLam}},\veg j_{\mat H},\veg j_{\mSig}, \veg m_{\mLam})$} \\
        \hline
        SEDR  &\makecell{ ({$\omega$}, $\omega$)\\({$\omega$}, \blue{$\omega$}) }&\makecell{ ({$\omega$}, $\omega$)\\({$\omega$}, \blue{$\omega$}) }&\makecell{ ($\omega^2$, $\omega^2$)\\({$\omega$}, $\omega$) }&\makecell{ ({$\omega$}, $\omega$)\\({$\omega$}, \blue{$\omega$}) }&\makecell{ ({$\omega$}, $\omega$)\\({$\omega$}, \blue{$\omega$}) }&\makecell{ ({$\omega$}, $\omega$)\\({$\omega$}, \blue{$\omega$}) }&\makecell{$\Re(\veg j_{\mLam},\veg j_{\mat H},\veg m_{\mLam},\veg m_{\mat H},\veg m_{\mSig})$ \\$\Im( \veg j_{\mLam},\veg j_{\mat H},\veg j_{\mSig},\veg m_{\mLam},\veg m_{\mat H},\veg m_{\mSig})$}  \\

        \hline

        \multicolumn{8}{|c|}{\textbf{(e) Electric exterior near field}} \\ 
        \hline
        Regime & \makecell{$\veg E(\Re \veg j_{\mLam})$\\$\veg E(\Im \veg j_{\mLam})$} & \makecell{$\veg E(\Re \veg j_{\mat H})$\\$\veg E(\Im \veg j_{\mat H})$} & \makecell{$\veg E(\Re \veg j_{\mSig})$\\$\veg E(\Im \veg j_{\mSig})$} & \makecell{$\veg E(\Re \veg m_{\mLam})$\\$\veg E(\Im \veg m_{\mLam})$} & \makecell{$\veg E(\Re \veg m_{\mat H})$\\$\veg E(\Im \veg m_{\mat H})$} & \makecell{$\veg E(\Re \veg m_{\mSig})$\\$\veg E(\Im \veg m_{\mSig})$} & Current dominant components \\
        \hline
        QSR  &\makecell{ ({$\omega^5$}, $\omega^2$)\\({$\omega^2$}, \blue{$\omega^5$}) }&\makecell{ ({$\omega^5$}, $\omega^2$)\\({$\omega^2$}, \blue{$\omega^5$}) }&\makecell{ ($\omega^3$, $1$)\\({$1$}, $\omega^3$) }&\makecell{ ({$1$}, $\omega^3$)\\({$\omega^3$}, \blue{$1$}) }&\makecell{ ({$\omega^2$}, $\omega^5$)\\({$\omega^5$}, \blue{$\omega^2$}) }&\makecell{ ({$\omega^2$}, $\omega^5$)\\({$\omega^5$}, \blue{$\omega^2$}) }&\makecell{$\Re(\veg j_{\mSig}, \veg m_{\mLam})$ \\$\Im(\veg j_{\mSig}, \veg m_{\mLam})$}\\
        \hline
        ECFR  &\makecell{ ({$\omega^{11/2}$}, $\omega^{5/2}$)\\({$\omega^2$}, \blue{$\omega^5$}) }&\makecell{ ({$\omega^6$}, $\omega^3$)\\({$\omega^2$}, \blue{$\omega^5$}) }&\makecell{ ($\omega^4$, $\omega$)\\({$1$}, $\omega^3$) }&\makecell{ ({$\omega^2$}, $\omega^5$)\\({$\omega^4$}, \blue{$\omega$}) }&\makecell{ ({$\omega^2$}, $\omega^5$)\\({$\omega^{11/2}$}, \blue{$\omega^{5/2}$}) }&\makecell{ ({$\omega^2$}, $\omega^5$)\\({$\omega^{11/2}$}, \blue{$\omega^{5/2}$}) }&$\Im(\veg j_{\mSig})$\\
        \hline
        SEDR  &\makecell{ ({$\omega^{5}$}, $\omega^{2}$)\\({$\omega^2$}, \blue{$\omega^5$}) }&\makecell{ ({$\omega^5$}, $\omega^2$)\\({$\omega^2$}, \blue{$\omega^5$}) }&\makecell{ ($\omega^4$, $\omega$)\\({$1$}, $\omega^3$) }&\makecell{ ({$\omega^2$}, $\omega^5$)\\({$\omega^5$}, \blue{$\omega^2$}) }&\makecell{ ({$\omega^2$}, $\omega^5$)\\({$\omega^{5}$}, \blue{$\omega^{2}$}) }&\makecell{ ({$\omega^2$}, $\omega^5$)\\({$\omega^{5}$}, \blue{$\omega^{2}$}) }&$\Im(\veg j_{\mSig})$\\
        \hline

        \multicolumn{8}{|c|}{\textbf{(f) Magnetic exterior near field}} \\ 
        \hline
        Regime & \makecell{$\veg H(\Re \veg j_{\mLam})$\\$\veg H(\Im \veg j_{\mLam})$} & \makecell{$\veg H(\Re \veg j_{\mat H})$\\$\veg H(\Im \veg j_{\mat H})$} & \makecell{$\veg H(\Re \veg j_{\mSig})$\\$\veg H(\Im \veg j_{\mSig})$} & \makecell{$\veg H(\Re \veg m_{\mLam})$\\$\veg H(\Im \veg m_{\mLam})$} & \makecell{$\veg H(\Re \veg m_{\mat H})$\\$\veg H(\Im \veg m_{\mat H})$} & \makecell{$\veg H(\Re \veg m_{\mSig})$\\$\veg H(\Im \veg m_{\mSig})$} & Current dominant components \\
        \hline
        QSR  &\makecell{ ({$\omega$}, $\omega^{4}$)\\({$\omega^4$}, \blue{$\omega$}) }&\makecell{ ({$\omega$}, $\omega^4$)\\({$\omega^4$}, \blue{$\omega$}) }&\makecell{ ($\omega$, $\omega^4$)\\({$\omega^4$}, $\omega$) }&\makecell{ ({$\omega^4$}, $\omega$)\\({$\omega$}, \blue{$\omega^4$}) }&\makecell{ ({$\omega^6$}, $\omega^3$)\\({$\omega^{3}$}, \blue{$\omega^{6}$}) }&\makecell{ ({$\omega^4$}, $\omega$)\\({$\omega$}, \blue{$\omega^{4}$}) }&\makecell{$\Re({\veg j_{\mLam}}, \veg j_{\mat H},\veg j_{\mSig}, \veg m_{\mLam}, {\veg m_{\mSig}})$ \\$\Im({\veg j_{\mLam}}, \veg j_{\mat H},\veg j_{\mSig}, \veg m_{\mLam}, {\veg m_{\mSig}})$}\\
        \hline
        ECFR  &\makecell{ ({$\omega^{3/2}$}, $\omega^{9/2}$)\\({$\omega^4$}, \blue{$\omega$}) }&\makecell{ ({$\omega^2$}, $\omega^5$)\\({$\omega^4$}, \blue{$\omega$}) }&\makecell{ ($\omega^2$, $\omega^5$)\\({$\omega^4$}, $\omega$) }&\makecell{ ({$\omega^6$}, $\omega^3$)\\({$\omega^2$}, \blue{$\omega^5$}) }&\makecell{ ({$\omega^6$}, $\omega^3$)\\({$\omega^{7/2}$}, \blue{$\omega^{13/2}$}) }&\makecell{ ({$\omega^4$}, $\omega$)\\({$\omega^{3/2}$}, \blue{$\omega^{9/2}$}) }&\makecell{$\Re( \veg m_{\mSig})$ \\$\Im({\veg j_{\mLam}},\veg j_{\mat H},\veg j_{\mSig})$}  \\
        \hline
        SEDR  &\makecell{ ({$\omega$}, $\omega^{4}$)\\({$\omega^4$}, \blue{$\omega$}) }&\makecell{ ({$\omega$}, $\omega^4$)\\({$\omega^4$}, \blue{$\omega$}) }&\makecell{ ($\omega^2$, $\omega^5$)\\({$\omega^4$}, $\omega$) }&\makecell{ ({$\omega^6$}, $\omega^3$)\\({$\omega^3$}, \blue{$\omega^6$}) }&\makecell{ ({$\omega^6$}, $\omega^3$)\\({$\omega^{3}$}, \blue{$\omega^{6}$}) }&\makecell{ ({$\omega^4$}, $\omega$)\\({$\omega$}, \blue{$\omega^{4}$}) }&\makecell{$\Re(\veg j_{\mLam},\veg j_{\mat H}, \veg m_{\mSig})$ \\$\Im({\veg j_{\mLam}},\veg j_{\mat H},\veg j_{\mSig},\veg m_{\mSig})$}  \\
        \hline
        
        \multicolumn{8}{|c|}{\textbf{(g) Rescaled current density}} \\ 
        \hline
        Regime & $b_R^{-1} \veg j_{\mLam}$ & $a_R^{-1} \veg j_{\mat H}$ & $a_R^{-1} \veg j_{\mSig}$ & $d_R^{-1} \veg m_{\mLam}$ & $c_R^{-1} \veg m_{\mat H}$ & $c_R^{-1} \veg m_{\mSig}$ & Recovered components \\
        \hline
        QSR & ($1,1$) & ($1,1$) & ($1,1$) & ($1,1$) & ($1,1$) & ($1,1$) & \makecell{$\Re( \veg j_{\mLam},\veg j_{\mat H},\veg j_{\mSig},\veg m_{\mLam},\veg m_{\mat H}, \veg m_{\mSig})$ \\$\Im( \veg j_{\mLam},\veg j_{\mat H},\veg j_{\mSig},\veg m_{\mLam},\veg m_{\mat H}, \veg m_{\mSig})$} \\
        \hline
        ECFR & ($\omega^{2}, \omega^{3/2}$) & ($\omega^{5/2}, \omega^{3/2}$) & ($\omega^{5/2}, \omega^{3/2}$) & ($\omega^{5/2}, \omega^{3/2}$) & ($\omega^{3/2}, \omega^{2}$) & ($\omega^{3/2}, \omega^{2}$) & \makecell{$\Re( \veg j_{\mLam},\veg m_{\mat H},\veg m_{\mSig})$ \\$\Im(\veg j_{\mLam},\veg j_{\mat H},\veg j_{\mSig}, \veg m_{\mLam})$} \\
        \hline
        SEDR & ($\omega^{3/2}, \omega^{3/2}$) & ($\omega^{3/2}, \omega^{3/2}$) & ($\omega^{5/2}, \omega^{3/2}$) & ($\omega^{3/2}, \omega^{3/2}$) & ($\omega^{3/2}, \omega^{3/2}$) & ($\omega^{3/2}, \omega^{3/2}$) & \makecell{$\Re(\veg j_{\mLam},\veg j_{\mat H}, \veg m_{\mLam},\veg m_{\mat H},\veg m_{\mSig})$ \\$\Im( \veg j_{\mLam},\veg j_{\mat H},\veg j_{\mSig},\veg m_{\mLam},\veg m_{\mat H},\veg m_{\mSig})$} \\
        \hline
    \end{tabular}
\end{table*}

\section{A new Preconditioned PMCHWT}
\label{sec:precondForm}
This section aims at defining a preconditioned formulation of the PMCHWT equation with the following objectives: (i) stabilizing the condition number with respect to both the electrical length and the conductivity of the scatterer in the regimes analyzed above, and (ii) preventing the loss of dominant components of the current in the low-frequency limit, in order to guarantee the accuracy of the solution and fields. This can be obtained by properly rescaling the quasi-Helmholtz decomposed components of the system, by means of primal and dual quasi-Helmholtz projectors. We propose here a quasi-Helmholtz projectors based preconditioning strategy tailored to the magnetic frill (or loop) excitations analyzed in tables \ref{table:indtablecurrent} and \ref{table:captablecurrent}, referred to as excitations of inductive (ind.) or capacitive (cap.) types. The proposed formulation follows from a left-right preconditioning of the original one, i.e.,
\begin{equation}
    \mat L \mat Z \mat R \vec y = \mat{L} \, \begin{pmatrix}
    \vec e \\ \vec h
\end{pmatrix}\,,
\end{equation}
where the left and right preconditioning matrices take the form
\begin{equation}
    \mat L =
    \begin{cases}
        \begin{pmatrix}
        \frac{1}{\sqrt{\eta_0}}\mat G^{-\T} {\mathbb M}_{\text{up}} \mat G^{-1} & \boldsymbol{0} \\
        \boldsymbol{0} & \sqrt{\eta_0}\mat G^{-\T}{\mathbb M}_{\text{low}} \mat G^{-1}
        \end{pmatrix} & \text{ind. type} \\
        \begin{pmatrix}
        \frac{1}{\sqrt{\eta_0}}{\mat M}_{\text{up}} & \boldsymbol{0} \\
        \boldsymbol{0} & \sqrt{\eta_0}{\mat M}_{\text{low}}
        \end{pmatrix} & \text{cap. type}
    \end{cases}
    \label{eqn:matL}
\end{equation}
\begin{equation}
    \mat R =
    \begin{cases}
        \begin{pmatrix}
        \frac{1}{\sqrt{\eta_0}}{\mat M}_{\text{left}} & \boldsymbol{0} \\
        \boldsymbol{0} & \sqrt{\eta_0}{\mat M}_{\text{right}}
        \end{pmatrix} & \text{ind. type}\\
        \begin{pmatrix}
        \frac{1}{\sqrt{\eta_0}}\mathbb{G}^{-1}{\mathbb M}_{\text{left}} \mat \mathbb{G}^{-\T} & \boldsymbol{0} \\
        \boldsymbol{0} & \sqrt{\eta_0}\mathbb{G}^{-1}{\mathbb M}_{\text{right}} \mat \mathbb{G}^{-\T}
        \end{pmatrix} & \text{cap. type}
    \end{cases}
    \label{eqn:matR}
\end{equation}
As already noticed in \Cref{sec:dominantexcitation}, the frequency behavior of the harmonic component of the electric field produced by a magnetic frill is different in case the excitation is of inductive or capacitive type. To correctly enforce the required cancellations, we employ two different, and symmetric, preconditioning strategies for the inductive and capacitive cases (equations \eqref{eqn:matL}, \eqref{eqn:matR}).
Indeed, the use of dual projectors in the left preconditioning matrix and primal projectors in the right preconditioning matrix for inductive excitations, and vice versa for capacitive excitations, allows to correctly enforce the cancellation of the static part of the electric component of the RHS projected on the rotated non-solenoidal subspace, non including the quasi-harmonic subspace, if the excitation is of inductive type, and the cancellation of the static part of the electric component of the RHS projected on the solenoidal subspace, including the quasi-harmonic subspace, if the excitation is capacitive. 

The preconditioning blocks used for inductive type excitations are
\begin{align} 
{\mathbb{M}}_{\text{up}} &=  {a_L}\, \PSH +{b_L}\, \PL\,, \\ 
{\mathbb{M}}_{\text{low}} &= {c_L} \, \PSH +{d_L} \, \PL\,, \\ 
{\mat{M}}_{\text{left}} &= {a_R} \, \PLH+ {b_R}  \,\PS\,,\\
{\mat{M}}_{\text{right}} &={c_R} \,\PLH+ {d_R} \,\PS\,.
\end{align}
Symmetrically, the ones used in case of capacitive type excitations are
\begin{align} 
{\mat{M}}_{\text{up}} &= {a_L} \, \PLH+ {b_L}  \,\PS\,,\\
{\mat{M}}_{\text{low}} &={c_L} \,\PLH+ {d_L} \,\PS\,\\
{\mathbb{M}}_{\text{left}} &=  {a_R}\, \PSH +{b_R}\, \PL\,, \\ 
{\mathbb{M}}_{\text{right}} &= {c_R} \, \PSH +{d_R} \, \PL\,, \,.
\end{align}

The scalar coefficients $\{a_L\operatorname{-}d_R\}$ should be set as a function of the frequency-conductivity regime considered, in order to stabilize the formulation with respect to the frequency and the conductivity independently.

The coefficients to be used in case of inductive and capacitive excitations are reported in Table~\ref{table:coefficients}. In addition, arguments justifying these choices are delineated in the following.

\begin{table*}[!t]
\renewcommand{\arraystretch}{1.4}
\setlength\tabcolsep{5pt}
\caption{Scalar multiplicative coefficients defining the proposed preconditioning strategy.}
\label{table:coefficients}
\centering
\begin{tabular}{|c| c c | c c | c c |} 
 \hline
 & QSR-ind & QSR-cap & ECFR-ind & ECFR-cap & SEDR-ind & SEDR-cap \\
 \hline 
$a_L$ & 
$1$ & $\chi^{-2}$ & 
$\chi^{1/2}\gamma^{-1}$ & $\chi^{-1/2}$ & 
$\chi^{-1/2}$ & $\chi^{-1/2}$\\
$b_L$ & 
$\chi^{-2}$ & $1$ & 
$\chi^{-1/2}$ & $\chi^{3/2}$ & 
$\chi^{-1/2}$ & $\chi^{3/2}$\\
$c_L$ & 
$\chi^{-1}$ & $\chi^{-1}$ & 
$\chi^{1/2}$ & $\chi^{1/2}$ & 
$\gamma^{1/2}$ & $\gamma^{1/2}$\\
$d_L$ & 
$\chi^{-1}$ & $\chi^{-1}$ & 
$\chi^{-1/2}\gamma$ & $\chi^{1/2}$ & 
$\gamma^{1/2}$ & $\gamma^{1/2}$\\
$a_R$ & 
$\chi$ & $\chi$ & 
$\chi^{-1/2}$ & $\chi^{-1/2}$ & 
$\chi^{-1/2}$ & $\chi^{-1/2}$\\
$b_R$ & 
$\chi$ & $\chi$ & 
$\chi^{1/2}\gamma$ & $\chi^{-1/2}$ & 
$\chi^{3/2}$ & $\chi^{-1/2}$\\
$c_R$ & 
$1$ & $\chi^{2}$ & 
$\chi^{-1/2}\gamma$ & $\chi^{1/2}$ & 
$\gamma^{1/2}$ & $\gamma^{1/2}$\\
$d_R$ & 
$\chi^{2}$ & $1$ & 
$\chi^{1/2}$ & $\chi^{-3/2}\gamma^2$ & 
$\gamma^{1/2}$ & $\gamma^{1/2}$\\
 \hline
\end{tabular}
\end{table*}

The conditioning analysis of the preconditioned formulation will be based on the asymptotic scalings of the loop-star decomposition of the preconditioned matrix, defined as
\begin{equation}
    \mat Z_{P,\mat \Lambda \mat H \mat \Sigma} \coloneqq   
    \begin{pmatrix}
        \mat A^\T & \boldsymbol{0} \\ \boldsymbol{0} & \mat A^\T
    \end{pmatrix} \, \mat L \mat Z \mat R \,\begin{pmatrix}
        \mat A & \boldsymbol{0} \\ \boldsymbol{0} & \mat A
    \end{pmatrix}\,.
\end{equation}

\subsection{The Quasi-Static Regime}
\label{sec:QSR_precond}
In the QSR, the contribution of losses inside the material does not influence the overall scalings of the quasi-Helmholtz decomposed formulation. Indeed, the loop-star decomposition identified in \eqref{eqn:LSPMCHWT_QS} corresponds to the one of the PMCHWT equation applied to lossless dielectric media, already identified in \cite{beghein2017lowfrequency}.
\subsubsection{Inductive type excitation}
Different sets of coefficients allow for frequency and conductivity preconditioning of the formulation. One classical choice \cite{beghein2017lowfrequency} is for example
\begin{equation}
    a_L = c_L = b_R = d_R = {\chi}^{1/2}\,\,\text{and}\,\,
    b_L = d_L = a_R = c_R = {\chi}^{-1/2}\,,\nonumber
\end{equation}
leading to the Loop-Star decomposition
\begin{equation}
   \mat Z_{P,\mat \Lambda \mat H \mat \Sigma} = \mathcal{O} \begin{pmatrix}
 1 & 1 & \chi & \chi & \chi & 1 \\
 \chi & \chi & \chi^2 & \chi^2 & 1 & \chi\\
 \chi & \chi & 1 & 1 & 1 & \chi  \\
 \chi & \chi & 1 & 1 & 1 & \chi \\
 \chi^2 & 1 & \chi & \chi & \chi & \chi^2 \\
 1 & 1 & \chi & \chi & \chi & 1
    \end{pmatrix}
    %\label{eqn:}
\end{equation}
leading to a constant condition number behavior in the low-frequency limit \cite{beghein2017lowfrequency}.

However, for the kind of excitation considered here, this choice would not allow for the  preservation of all the relevant current components. Indeed, the components $\vec j_{\mLam}$, $\vec j_{\mat H}$, and $\vec m_{\mSig}$, both in their real and imaginary parts, which are useful to reconstruct the magnetic field inside and outside the scatterer, would be lost in numerical cancellation in the low-frequency limit.

To preserve the accuracy of all the required current components, we opt instead to the choice of coefficients reported in Table~\ref{table:coefficients}, resulting in the loop-star decomposition
\begin{equation}
   \mat Z_{P,\mat \Lambda \mat H \mat \Sigma} = \mathcal{O} \begin{pmatrix}
 1 & 1 & 1 & 1 & 1 & 1 \\
 \chi^2 & \chi^2 & \chi^2 & \chi^2 & 1 & \chi^2\\
 \chi^2 & \chi^2 & 1 & 1 & 1 & \chi^2  \\
 \chi^2 & \chi^2 & 1 & 1 & 1 & \chi^2 \\
 \chi^2 & 1 & 1 & 1 & 1 & \chi^2 \\
 1 & 1 & 1 & 1 & 1 & 1
    \end{pmatrix}\,,
    \label{eqn:precondQSRInd}
\end{equation}
indicating the favorable condition properties of the preconditioned matrix and that the matrix has a well defined low-frequency limit converging to the static problem. The preserved current components correspond to the ones needed for the correct fields reconstruction, as shown in table \ref{table:indtablecurrent}.

\subsubsection{Capacitive type excitation}
Similarly as above, we determine from table \ref{table:indtablecurrent} the components to be preserved and design the multiplicative preconditioner accordingly. Our choice of coefficients, in table \ref{table:coefficients}, results in 
\begin{equation}
   \mat Z_{P,\mat \Lambda \mat H \mat \Sigma} = \mathcal{O} \begin{pmatrix}
 1 & 1 & 1 & 1 & \chi^2 & 1 \\
 1 & 1 & 1 & 1 & 1 & 1\\
 \chi^2 & \chi^2 & 1 & 1 & \chi^2 & \chi^2  \\
 \chi^2 & \chi^2 & 1 & 1 & \chi^2 & \chi^2 \\
 \chi^2 & 1 & 1 & 1 & \chi^2 & \chi^2 \\
 1 & 1 & 1 & 1 & \chi^2 & 1
    \end{pmatrix}\,,
    %\label{eqn:}
\end{equation}
which indicates that the resulting matrix has favorable conditioning properties, and allows to preserve the required current components (see Table~\ref{table:captablecurrent}).

\subsection{The Eddy-Current-Free Eddy-Current Regime}
\label{sec:ECFR_precond}
\subsubsection{Inductive type excitation}
One first intuitive choice for the scalar coefficients $\{a_L\operatorname{-}d_R\}$ is
\begin{align}
    a_L = c_L = b_R = d_R = {\chi}^{1/2} &\,,\quad\quad
    b_L = a_R = {\chi}^{-1/2}\,,\nonumber\\ \text{and}\,\,d_L = c_R &= \chi^{-1/2}\gamma\,.\nonumber
\end{align}
This would allow to retrieve all the dominant current components.
Moreover, this would result in the decomposition
\begin{equation}
    \mat Z_{P,\mat \Lambda \mat H \mat \Sigma} = \mathcal{O} \begin{pmatrix}
 1 & 1 & \chi &  
 \xi &  \xi & 1 \\
 \chi & \chi & \chi^2 &  
 \chi\xi & \gamma & \chi\\
 \chi & \chi & 1 & 
 \gamma & \gamma & \chi  \\
 \xi & \xi & \gamma & 
 1 & 1 &  \xi \\
 \xi^2 & 1 & \chi &  
 \xi &  \xi & \xi^2 \\
 1 & 1 & \chi &  
 \xi &  \xi & 1 
    \end{pmatrix}\,,
    \label{eqn:loopstarECFR_nonok}
\end{equation}

In presence of multiply connected geometries with finite genus $g$ instead, the preconditioned matrices exhibit a finite-dimensional nullspace in the static limit, of size $2g$. This can be seen for example by studying the loop-star decomposition for a genus-$1$ geometry and separating the toroidal from the poloidal global loop contributions (see \cite{bogaert2011low} for further details on this technique), i.e., by replacing matrix $\mat A$ with $\mat A' \coloneqq \left( \mat\Lambda \quad \mat H^T \quad \mat H^P \quad \mat \Sigma \right)$, where $\mat H^T$ and $\mat H^P$ are the transformation matrices from the toroidal and poloidal global loops to the RWG subspace. The decomposition of the preconditioned matrix reads
\begin{align}
    &\begin{pmatrix}
        \mat A'^\T & \boldsymbol{0} \\ \boldsymbol{0} & \mat A'^\T
    \end{pmatrix} \, \mat L \mat Z \mat R \,\begin{pmatrix}
        \mat A' & \boldsymbol{0} \\ \boldsymbol{0} & \mat A'
    \end{pmatrix} =\nonumber\\ &\quad\mathcal{O} \begin{pmatrix}
 1 & 1 &1 & \chi &  
 \xi &  \xi&  \xi & 1 \\
 \chi & \chi& \chi & \chi^2 &  
 \chi\xi & \boldsymbol{\xi^2\gamma}& \gamma & \chi\\
 \chi & \chi& \chi & \chi^2 & 
 \chi\xi & \boldsymbol{\xi^2\gamma}& \boldsymbol{\xi^2\gamma} & \chi\\
 \chi & \chi& \chi & 1 & 
 \gamma & \gamma& \gamma & \chi\\
 \xi & \xi& \xi & \gamma & 
 1 & 1 & 1 &  \xi \\
 \xi^2 & \boldsymbol{\xi^2} & 1 & \chi &  
 \xi &  \xi&  \xi & \xi^2 \\
 \xi^2 & \boldsymbol{\xi^2} & \boldsymbol{\xi^2}& \chi &  
 \xi &  \xi &\xi & \xi^2 \\
 1 & 1 & 1& \chi &  
 \xi &  \xi &\xi & 1
    \end{pmatrix}\,,
    \label{eqn:ecfr_bogaert_analysis}
\end{align}
where the terms in bold take into account the fact that the static part of the MFIO matrix vanishes when the basis loop does not pass through the test loop \cite{bogaert2011low}.
The application of \eqref{eqn:ecfr_bogaert_analysis} to the $V$-dimensional (with $V$ denoting the number of vertices of the mesh) space of magnetic currents composed of local RWG loops and the toroidal global loop, that is, the vector
\begin{equation}
    \begin{pmatrix}
        \vec 0 & \vec 0 & \vec 0 & \vec 0 & \vec m_{\mat \Lambda} & \vec m_{\mat H^T} & \vec 0 & \vec 0
    \end{pmatrix}^\T
\end{equation}
results in
\begin{equation}
    \begin{pmatrix}
        \mathcal{O}(\xi) & \mathcal{O}(\chi\xi) & \mathcal{O}(\chi\xi) & \mathcal{O}(\gamma) & \vec v & \mathcal{O}(\xi) & \mathcal{O}(\xi) & \mathcal{O}(\xi)
    \end{pmatrix}^\T\,,
\end{equation}
where $\vec v$ is given by the product between a rank-deficient $(V-1)\times V$ block and the $V$-dimensional current. As $\chi \ll 1$, $\gamma \ll 1$, and $\xi \ll 1$ in the ECFR, the nullspace of such rank-deficient block corresponds to the nullspace of the entire matrix in the low-frequency limit.
Similar results hold for the $V-$dimensional magnetic current space composed of local RWG loops and the poloidal global loop, resulting, overall, in a two-dimensional nullspace of the preconditioned formulation, which, numerically, takes the form of two isolated singular values decaying toward lower frequencies.
Note that a similar analysis applied to the preconditioned formulation in the QSR confirms the absence of a nullspace.

To remove the artificial nullspace introduced by this suboptimal preconditioning, we further manipulate the set of coefficients $\{a_L\operatorname{-}d_R\}$. The choice reported in table \ref{table:coefficients} on the one hand allows the retrieval of all the useful current components---as shown in table \ref{table:indtablecurrent}---and on the other hand results in the behavior
\begin{equation}
    \mat Z_{P,\mat \Lambda \mat H \mat \Sigma} = \mathcal{O} \begin{pmatrix}
 1 & 1 & \chi\gamma &  
 \xi &  \xi & 1 \\
 \xi & \xi & \chi^2 &  
 \xi^2 & 1 & \xi\\
 \xi & \xi & 1 & 
 1 & 1 & \xi  \\
 \xi & \xi & \gamma^2 & 
 1 & 1 &  \xi \\
 \xi^2 & 1 & \chi\gamma &  
 \xi &  \xi & \xi^2 \\
 1 & 1 & \chi\gamma &  
 \xi &  \xi & 1 
    \end{pmatrix}\,,
    \label{eqn:loopstarEDFR_ok}
\end{equation}
indicating favorable conditioning properties.
Moreover, by applying an analysis similar to the one above, 
the absence of a nullspace in the low-frequency limit related to the non-zero genus of the geometry can be verified.

\subsubsection{Capacitive type excitation}
An effective choice of coefficients to retrieve all the useful current components is the one reported in table \ref{table:coefficients}, resulting in the loop-star decomposition
\begin{equation}
    \mat Z_{P,\mat \Lambda \mat H \mat \Sigma} = \mathcal{O} \begin{pmatrix}
 1 & 1 & 1 &  
 1 &  \xi^2 & 1 \\
 1 & 1 & 1 &  
 1 & 1 & 1\\
 \chi^2 & \chi^2 & 1 & 
 \gamma^2 & \chi^2 & \chi^2  \\
 \xi^2 & \xi^2 & 1 & 
 1 & \xi^2 &  \xi^2 \\
 \xi^2 & 1 & 1 &  
 1 &  \xi^2 & \xi^2 \\
 1 & 1 & 1 &  
 1 &  \xi^2 & 1 
    \end{pmatrix}\,,
    %\label{eqn:}
\end{equation}
indicating that the formulation has a well-defined low-frequency limit and favorable conditioning properties.

\subsection{The Skin-Effect-Dominated Eddy-Current Regime}
\subsubsection{Inductive type excitation}
\label{sec:SEDR_precond}
As indicated in Table~\ref{table:indtablecurrent}, all the current components, in real and imaginary parts, are required to correctly retrieve the fields in the low-frequency limit. This requirement is satisfied by the choice of coefficients in Table~\ref{table:coefficients}, resulting in the following behavior of the loop-star decomposed system
\begin{align}
    &\mat Z_{P,\mat \Lambda \mat H \mat \Sigma}=\nonumber\\&\mathcal{O} \begin{pmatrix}
 1 & 1 & \chi^2 &  
 \xi^{-1/2} &  \xi^{-1/2} & \xi^{-1/2} \\
 1 & 1 & \chi^2 &  
 \xi^{-1/2} &  \xi^{-1/2} & \xi^{-1/2}\\
 1 & 1 & 1 & 
 \xi^{-1/2} &  \xi^{-1/2} & \xi^{-1/2}  \\
 \xi^{-1/2} & \xi^{-1/2} & \chi^{3/2}\gamma^{1/2} & 
 1 & 1 &  1 \\
 \xi^{-1/2} & \xi^{-1/2} & \chi^{3/2}\gamma^{1/2} &  
 1 &  1 & 1 \\
 \xi^{-1/2} & \xi^{-1/2} & \chi^{3/2}\gamma^{1/2} &  
 1 &  1 & 1 
    \end{pmatrix}\,,
    \label{eqn:prec_SEDR}
\end{align}
attesting to the favorable conditioning properties of the formulation.

\subsubsection{Capacitive type excitation}
In this case, the preconditioning strategy implemented for the inductive type excitation can be symmetrized, by exchanging the role of left and right coefficients (see Table~\ref{table:coefficients}), resulting in the same preconditioned formulation as in \eqref{eqn:prec_SEDR}, 
which also allows for the retrieval of the required current components.

\section{Numerical Results}
\label{sec:num_res}

\begin{figure*}
\subfloat[\label{}]{%
  \includegraphics[width=1
\columnwidth]{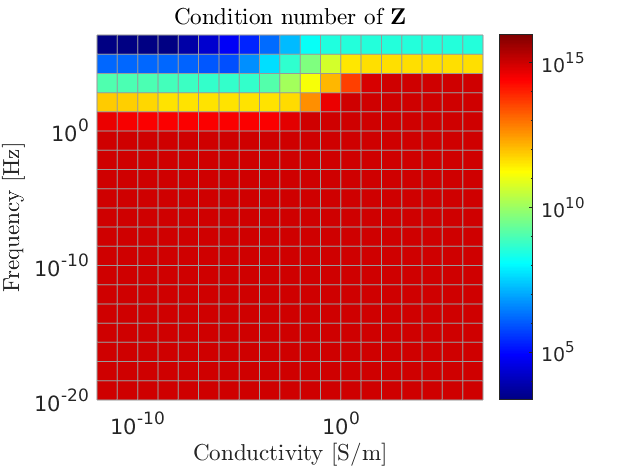}
}
\hfill
\subfloat[\label{}]{%
\includegraphics[width=1\columnwidth]{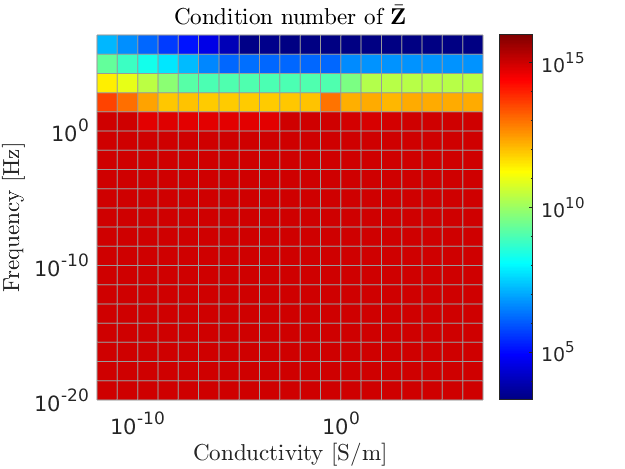}
}
\hfill
\subfloat[\label{}]{%
\includegraphics[width=1\columnwidth]{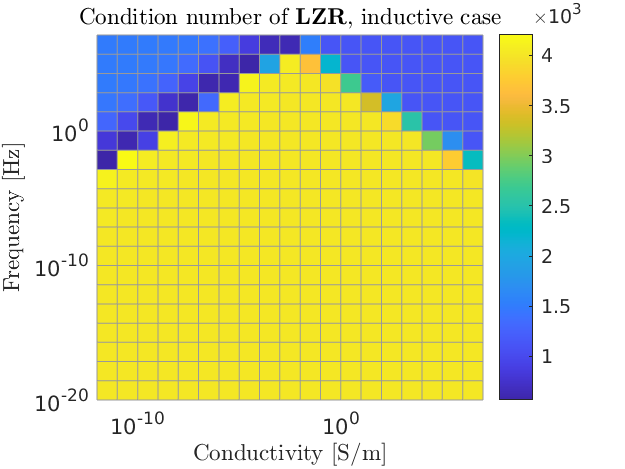}
}
\hfill
\subfloat[\label{}]{%
\includegraphics[width=1\columnwidth]{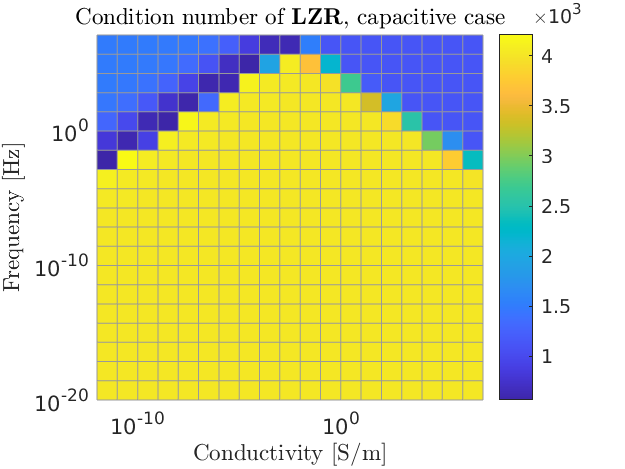}
}
\caption{Condition number of the formulations applied to torus with major and minor radii of \SI{1}{\meter} and \SI{0.3}{\meter} discretized at $h\simeq \SI{0.44}{\meter}$. The preconditioning applied is tailored to inductive (c) and capacitive (d) excitations.}
\label{fig:cnmaps}
\end{figure*}

In this section,  numerical results will illustrate the favorable spectral properties, the stability, and the accuracy of the proposed formulation.
First, the proposed formulation is applied to a multiply connected toroidal structure and its condition number is compared to that of the PMCHWT equation (matrix $\mat Z$) and of its rescaled version
(matrix $\bar{\mat Z}$)
%The condition number of the new formulation, matrix $\mat L \mat Z \mat R$,
in a wide range of conductivities and frequencies (always in the low-frequency regime), both for inductive and capacitive-type excitation versions (\Cref{fig:cnmaps}). On the one hand the low-frequency breakdown of the standard PMCHWT matrices $\mat Z$ or $\bar{\mat Z}$ is clear (as well as other conductivity related instabilities mentioned above), on the other hand the stability of the proposed formulation, that exhibits a condition number that remains constant in frequency and conductivity in the three regimes considered.

\begin{figure}
\centerline{\includegraphics[width=\columnwidth]{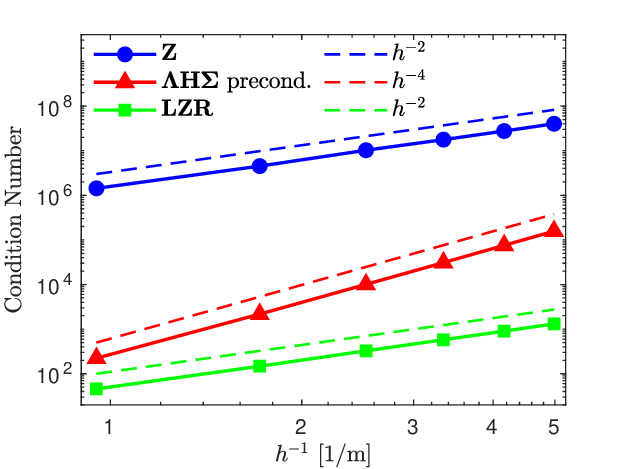}}
\caption{Condition number of the formulations applied to a sphere of radius \SI{1}{\meter} of conductivity \SI{10}{S/m} excited at \SI{1e5}{\hertz} at different mesh refinements: comparison between the original PMCHWT formulation ($\mat Z$), the Loop-Star preconditioned formulation ($\mat \Lambda \mat H \mat \Sigma$ precond.), and this work ($\mat L\mat Z\mat R$).}
\label{fig:cond_h}
\end{figure}

Next, we verify that the quasi-Helmholtz projectors based preconditioning does not degrade the spectral properties of the formulation with respect to mesh refinement by comparing the condition number of the standard PMCHWT matrix and of the proposed formulation with the one of a PMCHWT formulation preconditioned with a loop-star approach.
The rescaling coefficients employed in the loop-star preconditioning are the same as the one used in our formulation, in Table~\ref{table:coefficients}. The results, presented in Figure~\ref{fig:cond_h}, show that the condition number of matrix $\mat Z$ and the new formulation based on quasi-Helmholtz projector increase as $h^{-2}$ with refinement, which confirms the fact that the use of quasi-Helmholtz preconditioning is not detrimental with respect to the PMCHWT dense discretization breakdown. The condition number of the loop-star preconditioned formulation, instead, grows as $h^{-4}$.%testifying a detrimental impact of the application of $\mat \Lambda$ and $\mat \Sigma$ matrices, given their differential nature.

Beyond the conditioning properties of the equations under study, we next need to assess their accuracy. We consider two canonical circuital structures and employ our full-wave solver to evaluate their impedance, of both inductive and capacitive type. These values, derived from the voltages and currents themselves computed from the simulated total fields, are then compared to the values predicted from circuit theory to verify the accuracy of the schemes.

\begin{figure}
\centerline{\includegraphics[width=\columnwidth]{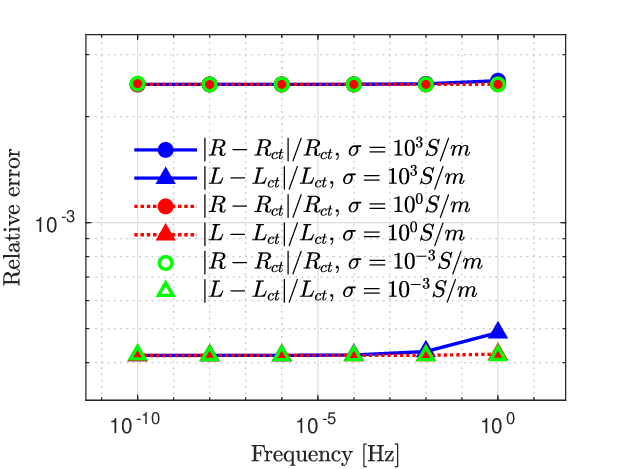}}
\caption{Full-wave evaluation of the impedance of a toroidal structure with major and minor radii of \SI{1}{\meter} and \SI{0.2}{\meter} of different conductivities: relative difference with respect to circuit theory expectations $R_{ct}$ and $L_{ct}$.}
\label{fig:ind_torus}
\end{figure}

\begin{figure}
\centerline{\includegraphics[width=\columnwidth]{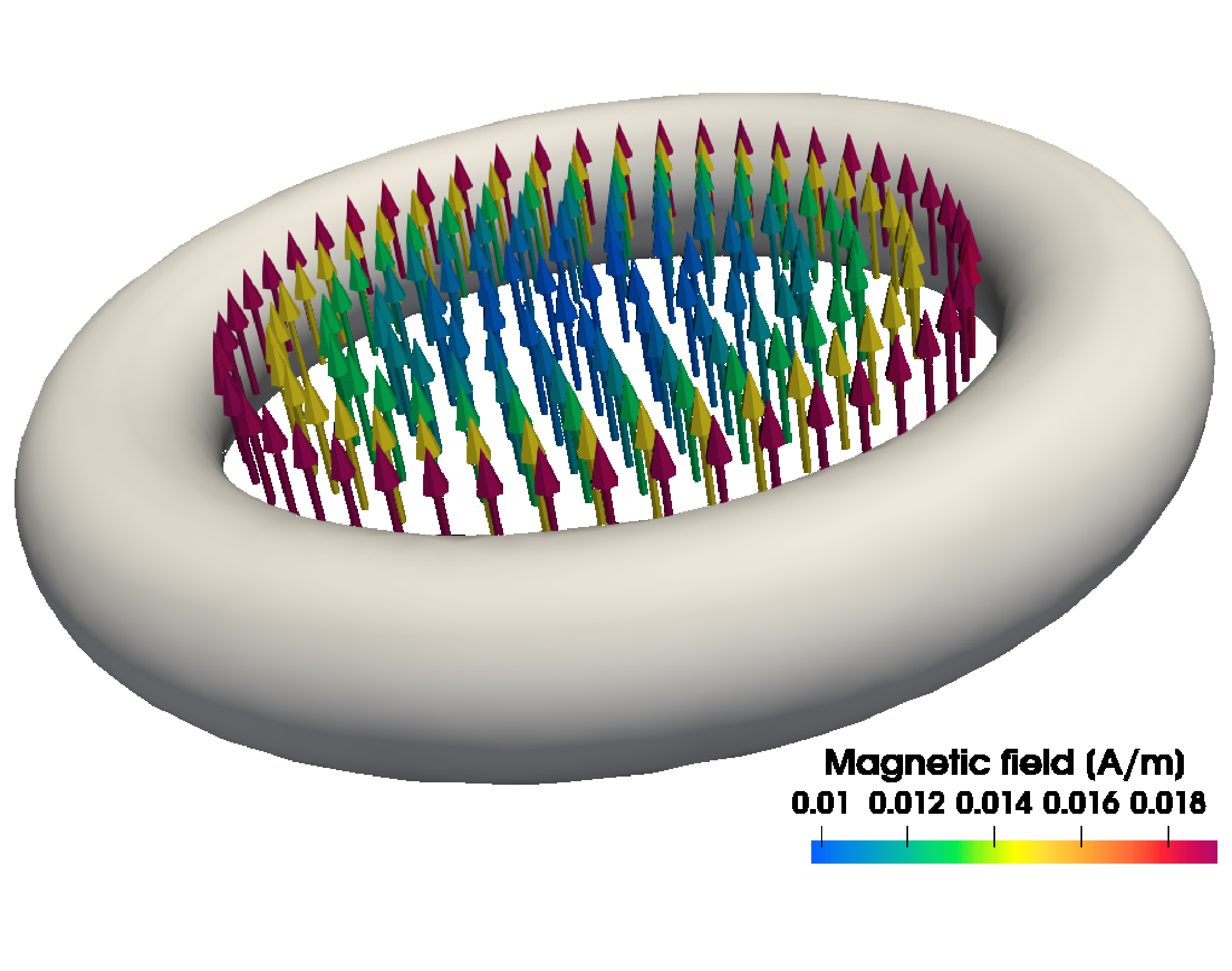}}
\caption{Torus of major and minor radii of \SI{1}{\meter} and \SI{0.2}{\meter} of conductivity \SI{1}{S/m} excited by a magnetic frill imposing \SI{1}{\volt} at \SI{10}{\milli\hertz}: induced magnetic field.}
\label{fig:torus_field}
\end{figure}

The first structure considered is a torus with major and minor radii $R_M=$ \SI{1}{\meter} and $R_m=$ \SI{0.2}{\meter} excited by an inductive magnetic frill passing through the global toroidal loop of the geometry, which imposes a potential difference of \SI{1}{\volt}.The conductivity of the toroidal body $\Omega$ varies from \SI{1e-3}{S/m} to \SI{1e3}{S/m}. We estimate the resistance and inductance of the circuit from the electromagnetic field scattered by the equivalent currents retrieved from the solution of our formulation and compare these values with the approximation valid in low frequency for conductors with uniform cross section $R_{ct}\coloneqq 2R_M/(R_m^2 \sigma)$ \cite{feynman2011feynman}, and $L_{ct} \simeq $ \SI{1.807}{\micro \henry} \cite[Eq. 4.27]{paul2010inductance}.
The relative difference between the results obtained from our full-wave formulation and the ones known from circuit theory, $R_{ct}$ and $L_{ct}$, is reported in \Cref{fig:ind_torus} for a wide range of frequencies. 
In particular, the relative error of the resistance is approximately \SI{2.5e-3}, the one of the inductance \SI{4.2e-4}, both of them independently of frequency and conductivity of the line. Figure~\ref{fig:torus_field} represents the magnetic field induced by the current circulating in the loop, parallel to the torus axis and symmetric with respect to the center.

\begin{figure}
\centerline{\includegraphics[width=\columnwidth]{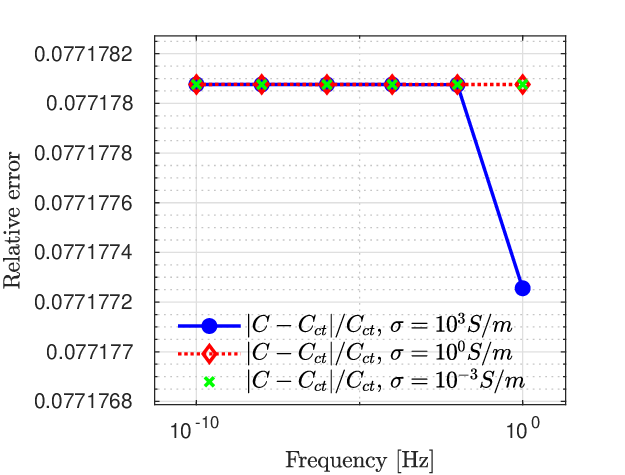}}
\caption{Full-wave evaluation of the impedance of a capacitive structure with circular parallel plates of radius \SI{4}{\meter} at distance \SI{0.2}{\meter} of different conductivities: relative difference with respect to circuit theory expectation $C_{ct}$.}
\label{fig:imp_cap}
\end{figure}

\begin{figure}
\centerline{\includegraphics[width=\columnwidth]{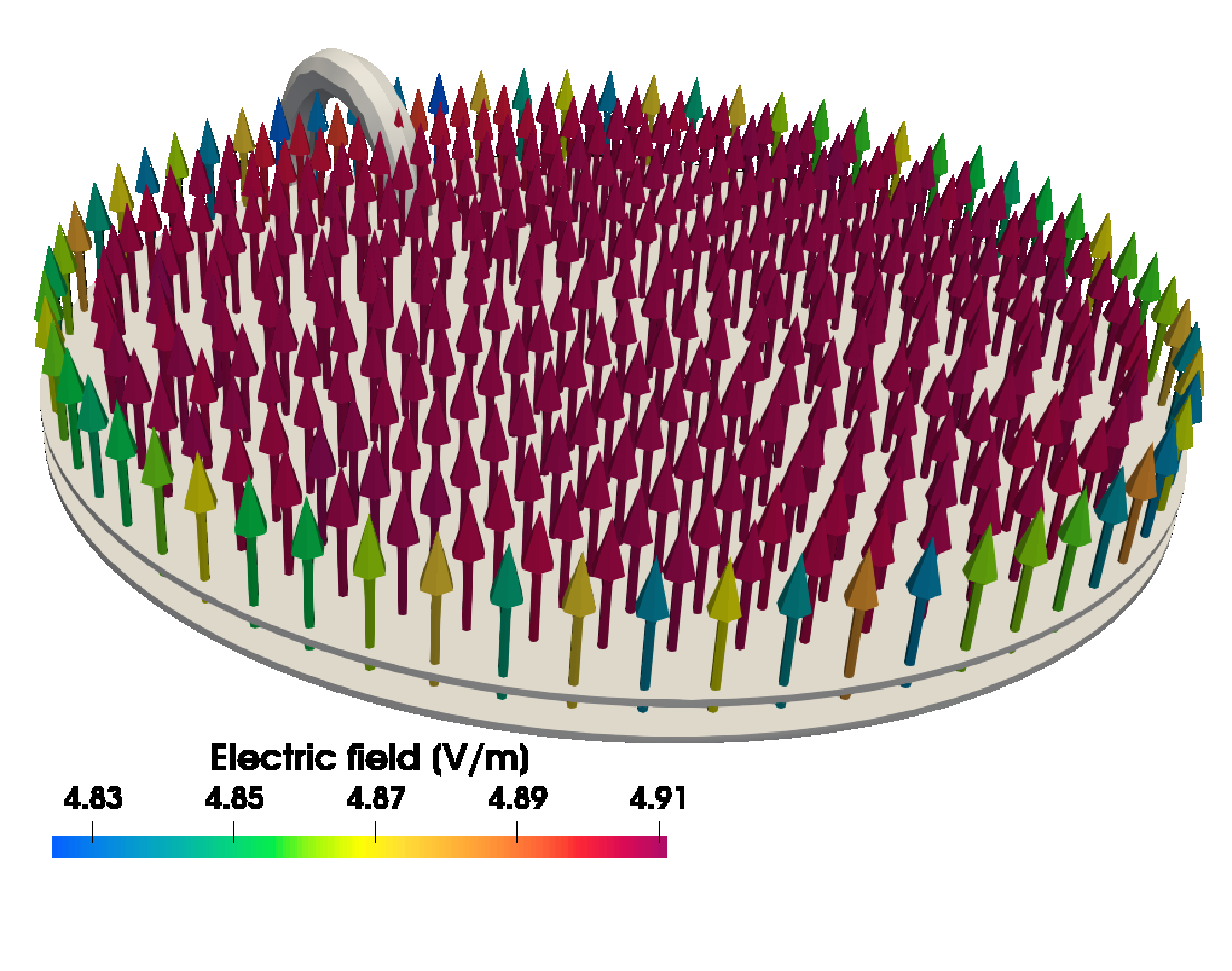}}
\caption{Capacitive structure with circular parallel plates of radius \SI{4}{\meter} at distance \SI{0.2}{\meter} of conductivity \SI{1e-3}{S/m} excited by a magnetic frill imposing \SI{1}{\volt} at \SI{1}{\hertz}: electric field between the plates.}
\label{fig:cap_field}
\end{figure}

The second structure instead is a parallel plates capacitor, made up of a couple of circular parallel plates of radius $R_p = \SI{4}{\meter}$ placed at distance $d = \SI{0.2}{\meter}$. The plates are connected by a wire of circular cross section of radius \SI{0.2}{\meter}. A potential difference of \SI{1}{\volt} is imposed across the plates using a capacitive magnetic frill around this wire. The capacitance of the structure is computed from the simulated fields, and compared to the value of capacitance expected from circuit theory $C_{ct} \coloneqq \epsilon_0\pi R_p^2/d$ in a wide range of frequencies and conductivities of the circuit. The results are summarized in \Cref{fig:imp_cap}, showing a relative difference between the two values of capacitance of approximately \SI{7.7e-2}, irrespectively of frequency and conductivity. Figure~\ref{fig:cap_field} shows a representation of the electric field generated between the plates, normally directed and uniform at a value near to \SI{5}{V/m}.

\begin{figure}
\subfloat[\label{}]{%
  \includegraphics[width=0.95
\columnwidth]{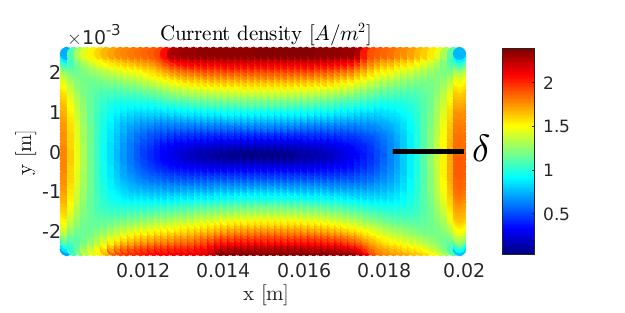}
}
\hfill
\subfloat[\label{}]{%
\includegraphics[width=0.95\columnwidth]{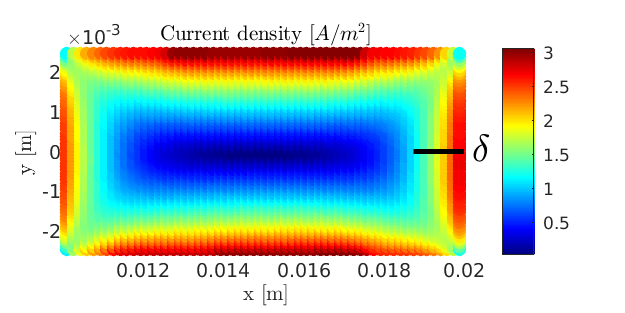}
}
\hfill
\subfloat[\label{}]{%
\includegraphics[width=0.95\columnwidth]{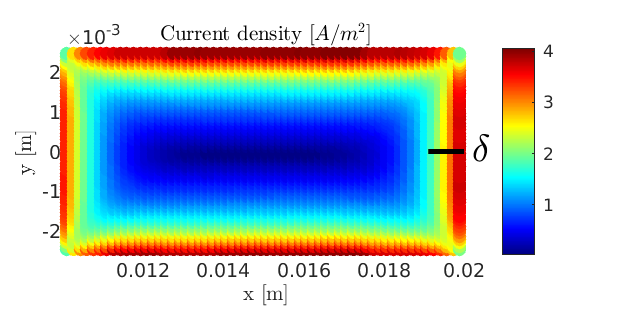}
}
\hfill
\subfloat[\label{}]{%
\includegraphics[width=0.95\columnwidth]{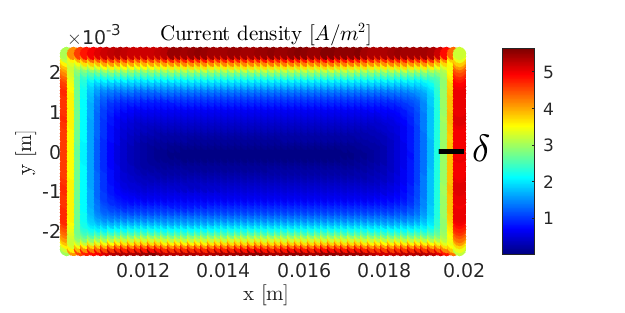}
}
\caption{Electric current density inside a rectangular cross section of sides \SI{1}{cm} and \SI{0.5}{cm} of conductivity \SI{1e7}{S/m} excited by an impinging plane wave travelling along $\hat{\veg y}$ at frequencies (a) \SI{8}{\kilo\hertz},(b) \SI{16}{\kilo\hertz},(c) \SI{32}{\kilo\hertz}, and (d) \SI{64}{\kilo\hertz}.}
\label{fig:skindepth}
\end{figure}

In the final numerical example, we demonstrate that our formulation can effectively track the field penetration depth variation inside a conductor. We considered in this case an highly conductive wire ($\sigma =  \SI{1e7}{S/m}$) of rectangular cross section enclosed in square shape with axis of symmetry directed toward $\hat{\veg y}$. After exciting the structure with a plane wave traveling along $\hat{\veg y}$, the resulting current density distribution in a cross section of the line at $z= \SI{0}{\meter}$ is shown in \Cref{fig:skindepth} for different frequencies of the impinging wave. We notice that the visible penetration length, which decreases when increasing the frequency, is comparable with the skin depth $\updelta$, which is also reported in the figures as the black lines for ease of comparison.

\section{Conclusion}
\label{sec:conclusion}
In this work we presented a novel full-wave formulation amenable to eddy current simulations in circuital frameworks. Built upon the PMCHWT integral equation, it results from a quasi-Helmholtz preconditioning strategy of projective type. With respect to the standard loop-star decomposition, the quasi-Helmholtz projectors offer several advantages, including their compatibility with fast solvers and acceleration strategies and their lack of detrimental effect on the condition number with respect to dense discretizations. The proposed preconditioning strategy offers the twofold advantage of curing the frequency and conductivity instabilities of the original formulation, including the severe low-frequency breakdown preventing its application toward statics, and curing the typical loss of accuracy occurring at very low-frequency. The resulting scheme is well-conditioned and stable for materials ranging from dielectric to highly conductive, applicable to both simply and multiply-connected structures and prone to seamless transitions between lower to higher frequencies, as well as to the modeling of multi-scale scenarios. Although the proposed scheme provides optimal accuracy in presence of magnetic frill type excitations, the procedure employed in defining the proposed preconditioning strategy can also be applied to other types of excitations.

\bibliographystyle{ieeetr}
\bibliography{private_vgiunzioni.bib}

\begin{thebibliography}{10}

\bibitem{balanis2012advanced}
C.~A. Balanis, {\em Advanced Engineering Electromagnetics}.
\newblock Hoboken, N.J: John Wiley \& Sons, 2nd ed~ed., 2012.

\bibitem{dirks1996quasistationary}
H.~K. Dirks, ``Quasi-stationary fields for microelectronic applications,'' {\em Electrical Engineering}, vol.~79, pp.~145--155, Apr. 1996.

\bibitem{sauter2011boundary}
S.~Sauter and C.~Schwab, {\em Boundary Element Methods}, vol.~39 of {\em Springer {{Series}} in {{Computational Mathematics}}}.
\newblock Berlin, Heidelberg: Springer Berlin Heidelberg, 2011.

\bibitem{steinbach2008numerical}
O.~Steinbach, {\em Numerical Approximation Methods for Elliptic Boundary Value Problems: Finite and Boundary Elements}.
\newblock New York: Springer, 2008.

\bibitem{poggio1973integral}
A.~Poggio and E.~Miller, ``Integral equation solutions of three-dimensional scattering problems,'' in {\em Computer {{Techniques}} for {{Electromagnetics}}}, pp.~159--264, Elsevier, 1973.

\bibitem{muller1958grundprobleme}
C.~Muller, ``Grundprobleme der mathematischen theorie elektromagnetischer schwingungen,'' {\em The American Mathematical Monthly}, vol.~65, p.~459, June 1958.

\bibitem{bogaert2014lowfrequency}
I.~Bogaert, K.~Cools, F.~Andriulli, and H.~Bagci, ``Low-frequency scaling of the standard and mixed magnetic field and {{M{\"u}ller}} integral equations,'' {\em IEEE Transactions on Antennas and Propagation}, vol.~62, pp.~822--831, Feb. 2014.

\bibitem{yla-oijala2008analysis}
P.~{Yl{\"a}-Oijala}, M.~Taskinen, and S.~J{\"a}rvenp{\"a}{\"a}, ``Analysis of surface integral equations in electromagnetic scattering and radiation problems,'' {\em Engineering Analysis with Boundary Elements}, vol.~32, pp.~196--209, Mar. 2008.

\bibitem{vantwout2022boundary}
E.~{van 't Wout}, S.~Haqshenas, P.~G{\'e}lat, T.~Betcke, and N.~Saffari, ``Boundary integral formulations for acoustic modelling of high-contrast media,'' {\em Computers \& Mathematics with Applications}, vol.~105, pp.~136--149, Jan. 2022.

\bibitem{adrian2021electromagnetic}
S.~Adrian, A.~Dely, D.~Consoli, A.~Merlini, and F.~Andriulli, ``Electromagnetic integral equations: Insights in conditioning and preconditioning,'' {\em IEEE Open Journal of Antennas and Propagation}, vol.~2, pp.~1143--1174, 2021.

\bibitem{cools2011Calderon}
K.~Cools, F.~Andriulli, and E.~Michielssen, ``A {{Calder{\'o}n}} multiplicative preconditioner for the {{PMCHWT}} integral equation,'' {\em IEEE Transactions on Antennas and Propagation}, vol.~59, pp.~4579--4587, Dec. 2011.

\bibitem{beghein2017lowfrequency}
Y.~Beghein, R.~Mitharwal, K.~Cools, and F.~Andriulli, ``On a low-frequency and refinement stable {{PMCHWT}} integral equation leveraging the quasi-{{Helmholtz}} projectors,'' {\em IEEE Transactions on Antennas and Propagation}, vol.~65, pp.~5365--5375, Oct. 2017.

\bibitem{hiptmair2007boundary}
R.~Hiptmair, ``Boundary element methods for eddy current computation,'' pp.~213--248, 2007.

\bibitem{kriezis1992eddy}
E.~Kriezis, T.~Tsiboukis, S.~Panas, and J.~Tegopoulos, ``Eddy currents: Theory and applications,'' {\em Proceedings of the IEEE}, vol.~80, no.~10, pp.~1559--1589, Oct./1992.

\bibitem{garcia-martin2011nondestructive}
J.~{Garc{\'i}a-Mart{\'i}n}, J.~{G{\'o}mez-Gil}, and E.~{V{\'a}zquez-S{\'a}nchez}, ``Non-destructive techniques based on eddy current testing,'' {\em Sensors}, vol.~11, pp.~2525--2565, Feb. 2011.

\bibitem{rucker1995various}
W.~Rucker, R.~Hoschek, and K.~Richter, ``Various {{BEM}} formulations for calculating eddy currents in terms of field variables,'' {\em IEEE Transactions on Magnetics}, vol.~31, pp.~1336--1341, May 1995.

\bibitem{niino2012Calderon}
K.~Niino and N.~Nishimura, ``Calder{\'o}n preconditioning approaches for {{PMCHWT}} formulations for {{Maxwell}}'s equations,'' {\em International Journal of Numerical Modelling: Electronic Networks, Devices and Fields}, vol.~25, pp.~558--572, Sept. 2012.

\bibitem{suyan2010comparative}
{Su Yan}, {Jian-Ming Jin}, and {Zaiping Nie}, ``A comparative study of {{Calder{\'o}n}} preconditioners for {{PMCHWT}} equations,'' {\em IEEE Transactions on Antennas and Propagation}, vol.~58, pp.~2375--2383, July 2010.

\bibitem{chhim2020eddy}
T.~Chhim, A.~Merlini, L.~Rahmouni, J.~Ortiz~G., and F.~Andriulli, ``Eddy current modeling in multiply connected regions via a full-wave solver based on the quasi-{{Helmholtz}} projectors,'' {\em IEEE Open Journal of Antennas and Propagation}, vol.~1, pp.~534--548, 2020.

\bibitem{andriulli2012loopstar}
F.~P. Andriulli, ``Loop-star and loop-tree decompositions: Analysis and efficient algorithms,'' {\em IEEE Transactions on Antennas and Propagation}, vol.~60, pp.~2347--2356, May 2012.

\bibitem{andriulli2013wellconditioned}
F.~Andriulli, K.~Cools, I.~Bogaert, and E.~Michielssen, ``On a well-conditioned electric field integral operator for multiply connected geometries,'' {\em IEEE Transactions on Antennas and Propagation}, vol.~61, pp.~2077--2087, Apr. 2013.

\bibitem{hackbuschsparse}
W.~Hackbusch, ``A sparse matrix arithmetic based on {{$\mathcal{H}$-matrices}}. part i: Introduction to {{$\mathcal{H}$-matrices}},''

\bibitem{hackbusch2000sparse}
W.~Hackbusch and B.~N. Khoromskij, ``A sparse {{$\mathcal{H}$-matrix}} arithmetic.: Part {{II}}: Application to multi-dimensional problems,'' {\em Computing}, vol.~64, pp.~21--47, Feb. 2000.

\bibitem{bogaert2011low}
I.~Bogaert, K.~Cools, F.~P. Andriulli, and D.~De~Zutter, ``Low frequency scaling of the mixed {{MFIE}} for scatterers with a non-simply connected surface,'' in {\em 2011 {{International Conference}} on {{Electromagnetics}} in {{Advanced Applications}}}, (Torino, Italy), pp.~951--954, IEEE, Sept. 2011.

\bibitem{zhiguoqian2007generalized}
{Zhi Guo Qian}, {Weng Cho Chew}, and R.~Suaya, ``Generalized impedance boundary condition for conductor modeling in surface integral equation,'' {\em IEEE Transactions on Microwave Theory and Techniques}, vol.~55, pp.~2354--2364, Nov. 2007.

\bibitem{quarteroni2008numerical}
A.~Quarteroni and A.~Valli, {\em Numerical Approximation of Partial Differential Equations}, vol.~23 of {\em Springer {{Series}} in {{Computational Mathematics}}}.
\newblock Berlin, Heidelberg: Springer Berlin Heidelberg, 2008.

\bibitem{henderson1981deriving}
H.~V. Henderson and S.~R. Searle, ``On deriving the inverse of a sum of matrices,'' {\em SIAM Review}, vol.~23, pp.~53--60, Jan. 1981.

\bibitem{tsai1972numerical}
L.~Tsai, ``A numerical solution for the near and far fields of an annular ring of magnetic current,'' {\em IEEE Transactions on Antennas and Propagation}, vol.~20, pp.~569--576, Sept. 1972.

\bibitem{popovic1982analysis}
B.~Popovi{\'c}, M.~B. Dragovi{\'c}, and A.~R. Djordjevi{\'c}, {\em Analysis and Synthesis of Wire-Antennas}.
\newblock No.~2 in Electronic \& Electrical Engineering Research Studies, Chichester, Sussex ; New York: Research Studies Press, 1982.

\bibitem{feynman2011feynman}
R.~P. Feynman, {\em The {{Feynman}} Lectures on Physics. {{Volume}} 1: {{Mainly}} Mechanics, Radiation, and Heat}.
\newblock New York: Basic Books, the new millennium edition, paperback first published~ed., 2011.

\bibitem{paul2010inductance}
C.~R. Paul, {\em Inductance: Loop and Partial}.
\newblock Hoboken, NJ: Wiley, 2010.

\end{thebibliography}

\end{document}